\documentclass[12pt,reqno]{amsart}
\usepackage{amsmath,amsfonts,amsthm,amsopn,amssymb,extarrows}
\usepackage{cite,marginnote}
\usepackage{color,enumitem,graphicx}
\usepackage[colorlinks=true,urlcolor=blue,
citecolor=red,linkcolor=blue,linktocpage,pdfpagelabels,
bookmarksnumbered,bookmarksopen]{hyperref}
\usepackage[english]{babel}

\usepackage[left=2.9cm,right=2.9cm,top=2.8cm,bottom=2.8cm]{geometry}
\usepackage[hyperpageref]{backref}

\numberwithin{equation}{section}

\makeindex
\makeindex
\newtheorem{theorem}{Theorem}[section]
\newtheorem{definition}[theorem]{Definition}
\newtheorem{lemma}[theorem]{Lemma}
\newtheorem{corollary}[theorem]{Corollary}
\newtheorem{proposition}[theorem]{Proposition}
\newtheorem{remark}[theorem]{Remark}

\newcommand{\s}{\section}

\newcommand{\R}{\mathbb R}

\newcommand{\C}{\mathbb C}

\newcommand{\lab}{\label}
\newcommand{\bt}{\begin{theorem}}
	\newcommand{\et}{\end{theorem}}
\newcommand{\bl}{\begin{lemma}}
	\newcommand{\el}{\end{lemma}}
\newcommand{\bd}{\begin{definition}}
	\newcommand{\ed}{\end{definition}}
\newcommand{\bc}{\begin{corollary}}
	\newcommand{\ec}{\end{corollary}}
\newcommand{\bp}{\begin{proof}}
	\newcommand{\ep}{\end{proof}}
\newcommand{\bx}{\begin{example}}
	\newcommand{\ex}{\end{example}}
\newcommand{\bi}{\begin{exercise}}
	\newcommand{\ei}{\end{exercise}}
\newcommand{\bo}{\begin{proposition}}
	\newcommand{\eo}{\end{proposition}}
\newcommand{\br}{\begin{remark}}
	\newcommand{\er}{\end{remark}}
\newcommand{\beq}{\begin{equation}}
	\newcommand{\eeq}{\end{equation}}
\newcommand{\ba}{\begin{align}}
	\newcommand{\ea}{\end{align}}
\newcommand{\bn}{\begin{enumerate}}
	\newcommand{\en}{\end{enumerate}}
\newcommand{\bg}{\begin{aligned}}
	\newcommand{\bcs}{\begin{cases}}
		\newcommand{\ecs}{\end{cases}}

	\newcommand{\bean}{\begin{eqnarray*}}
		\newcommand{\eean}{\end{eqnarray*}}

	\def\C{\mathbb{C}}

	\def\R{\mathbb{R}}

	\def\bd{\mathrm{bd}\,}

	\newcommand{\cm}{{\mathcal M}}

	\newcommand{\cp}{{\mathcal P}}

	\title[Multiple normalized solutions]{Multiple normalized solutions for two coupled Gross-Pitaevskii equations with attractive interactions and mass constriants}

	\author[J.~J.~Zhang]{Jianjun Zhang}
	\author[X.~X.~Zhong]{Xuexiu Zhong}
	\author[J.~F.~Zhou]{Jinfang Zhou}

    \address[J.~J.~Zhang]{\newline\indent College of Mathematics and Statistics
        \newline\indent
        Chongqing Jiaotong University
        \newline\indent
        Xuefu, Nan'an, 400074, Chongqing, PR China}
    \email{\href{mailto:zhangjianjun09@tsinghua.org.cn}{zhangjianjun09@tsinghua.org.cn}}

	\address[X.~X.~Zhong]{\newline\indent South China Research Center for Applied Mathematics and Interdisciplinary Studies
		\newline\indent
		South China Normal University
		\newline\indent
		Guangzhou 510631, P. R. China}
	\email{\href{mailto:zhongxuexiu1989@163.com}{zhongxuexiu1989@163.com}}
	
	\address[J.~F.~Zhou]{\newline\indent School of Mathematical Sciences
		\newline\indent
		South China Normal University
		\newline\indent
		Guangzhou 510631, P. R. China}
	\email{\href{mailto:jinfangnlsqjdbm@hotmail.com}{jinfangnlsqjdbm@hotmail.com}}

\thanks{Xuexiu Zhong was supported by the NSFC (No.12271184), Guangdong Basic and Applied Basic Research Foundation (2021A1515010034),Guangzhou Basic and Applied Basic Research Foundation(2024A04J10001). Jianjun Zhang was supported by the NSFC (No.12371109)}

\begin{document}
		
\begin{abstract}
We are concerned with the following system of two coupled time-independent Gross-Pitaevskii equations
$$
\begin{cases}
-\Delta u+\lambda_1 u=\mu_1|u|^{p-2}u+\nu\alpha |u|^{\alpha-2}|v|^{\beta}u ~\hbox{in}~ \R^N,\\
-\Delta v+\lambda_2 v=\mu_2|v|^{q-2}v+\nu\beta |u|^{\alpha}|v|^{\beta-2}v ~\hbox{in}~ \R^N,
\end{cases}
$$
which arises in two-components Bose-Einstein condensates and involve attractive Sobolev subcritical or critical interactions, i. e., $\nu>0$ and $\alpha+\beta\leq 2^*$. This system is employed by seeking critical points of the associated variational functional with the constrained mass below
$$\int_{\mathbb{R}^N}|u|^2 {\rm d}x=a, \quad \int_{\mathbb{R}^N}|v|^2 {\rm d}x=b.$$ In the mass mixed case, i. e., $2<p<2+\frac{4}{N}<q<2^*$, for some suitable $a,b,\nu$ and $\beta$, the system above admits two positive solutions. In particular, in the case $\alpha+\beta<2^*$, using variational methods on the $L^2$-ball, two positive solutions are obtained, one of which is a local minimizer and the second one is a mountain pass solution.
\begin{flushleft}
\end{flushleft}
\begin{flushleft}
\end{flushleft}

\end{abstract}
\maketitle

\s{Introduction}
\subsection{Background and motivation} In this paper, we study the following two coupled time-independent Gross-Pitaevskii equations
\beq\label{eq:230905-1}
\begin{cases}
	-i\partial_t\Psi_1=\Delta \Psi_1+\mu_1|\Psi_1|^{p-2}\Psi_1+\nu\alpha|\Psi_1|^{\alpha-2}|\Psi_2|^{\beta}\Psi_1,\\
	-i\partial_t\Psi_2=\Delta \Psi_2+\mu_2|\Psi_2|^{q-2}\Psi_2+\nu\beta|\Psi_1|^{\alpha}|\Psi_2|^{\beta-2}\Psi_2,\\
	\Psi_j=\Psi_j (x,t)\in\C, (x,t)\in \R^N\times\R, j=1,2,
\end{cases}
\eeq
where $i$ is the imaginary unit, and $\mu_i, \nu>0$ represents the intraspecies interaction and interspecies interaction respectively. This system arises in various physical contexts, particularly in the mean-field modeling of binary mixtures such as Bose-Einstein condensates (BEC) or degenerate quantum gases (e.g., Bose-Fermi or Fermi-Fermi mixtures); see, for instance, \cite{Adhikari2007,Bagnato2015,Esry1997,Malomed2008}. A key physical feature of \eqref{eq:230905-1} is the conservation of masses, meaning that the $L^2$-norms: $$|\Psi_1(\cdot,t)|_2 ~\hbox{and}~ |\Psi_2(\cdot,t)|_2$$
are independent of $t\in\R$. These conserved quantities have distinct physical interpretations: they represent the number of particles for each component in Bose-Einstein condensates, or the beam power in the context of nonlinear optics.

The main purpose of this paper is focused on solitary wave solutions of \eqref{eq:230905-1}, i. e., $\Psi_1(x,t)=e^{i\lambda_1 t}u_1(x)$ and $\Psi_2(x,t)=e^{i\lambda_2 t}u_2(x)$. Substituting this into system \eqref{eq:230905-1} yields the elliptic system
\beq\label{eq:240722-1}
\begin{cases}
-\Delta u+\lambda_1 u=\mu_1|u|^{p-2}u+\nu\alpha|u|^{\alpha-2}u|v|^{\beta} ~\hbox{in}~\R^N,\\
-\Delta v+\lambda_2 v=\mu_2|v|^{q-2}v+\nu\beta|u|^{\alpha}|v|^{\beta-2}v ~\hbox{in}~\R^N.
\end{cases}
\eeq
A solution $(u,v)$ satisfying the prescribed mass constraints
\beq\label{eq:240722-2}
\|\Psi_1(\cdot,0)\|_2^2=\int_{\R^N}|u|^2{\rm d}x=a, \quad \quad \|\Psi_2(\cdot,0)\|_2^2=\int_{\R^N}|v|^2{\rm d}x=b,
\eeq
is called as a normalized solution in the literature. To find solutions of Eqs.\eqref{eq:240722-1}-\eqref{eq:240722-2}, a classical strategy is to seek critical points $(u,v)\in H^1(\R^N)\times H^1(\R^N)$ of the energy functional
$$I(u,v)=\frac{1}{2}\int_{\R^N}(|\nabla u|^2+|\nabla v|^2){\rm d}x-\frac{\mu_1}{p}\int_{\R^N}|u|^p{\rm d}x-\frac{\mu_2}{q}\int_{\R^N}|v|^q{\rm d}x-\nu\int_{\R^N}|u|^{\alpha}|v|^{\beta}{\rm d}x$$
under the constraint
$$
T(a,b):=\{(u,v)\in H^1(\R^N,\R^2): \|u\|_2^2=a,\|v\|_2^2=b\}.
$$
Here, $\lambda_1,\lambda_2$ appear as Lagrange multipliers associated with the constraints.

In recent years, normalized solutions to system \eqref{eq:240722-1}-\eqref{eq:240722-2} have received widespread attention. This interest stems from the fact that the presence of the $L^2$ constraint renders many methods developed for unconstrained variational problems (the \textit{fixed frequency} problem) inapplicable. Compared to the \textit{fixed frequency} problem, the study of normalized solutions introduces significant technical challenges within the variational framework. A primary difficulty is the lack of compactness for Palais-Smale sequences on the constraint manifold. This arises because the embedding $H^1(\R^N)\hookrightarrow L^2(\R^N)$ is non-compact. Even when restricting to the radial subspace $H_{rad}^{1}(\R^N)$, the embedding $H^{1}_{rad}(\R^N)\hookrightarrow L^2(\R^N)$ remains non-compact. For the \textit{fixed frequency} problem, a nontrivial weak limit is typically also a solution. In contrast, for the \textit{fixed mass} problem, even if the weak limit is nontrivial, it may fail to satisfy the constraint condition. Furthermore, it is well-known that the $L^2$-critical exponent $\bar p:=2+\frac{4}{N}$ plays a pivotal role in the study of normalized solutions, critically influencing the geometry of the associated functional.

If all exponents $p,q$ and $\alpha+\beta$ are smaller than $\bar p$, i.e., $2<p,q,\alpha+\beta<\bar p$, then system \eqref{eq:240722-1}-\eqref{eq:240722-2} is purely $L^2$-subcritical, and the energy functional $I(u,v)$ is bounded below on $T(a,b)$. In this case, Gou and Jeanjean \cite{Gou2016} established the existence of a global minimizer, which is also a ground state solution. Conversely, if at least one of the exponents $p,q,\alpha+\beta$ exceeds $\bar p$, the system becomes $L^2$-supercritical, and $I(u,v)$ is unbounded below on $T(a,b)$. For this supercritical case, relevant results can be found in works such as \cite{Bartsch2016,Bartsch2021,Bartsch2018,Jeanjean2024a,Gou2018,Bartsch2023,Guo2022,Li2021a,Vicentiu2023}. Specifically, within the $L^2-$supercritical case, if all exponents $p,q,\alpha+\beta$ are greater than $\bar p$, it is termed the purely $L^2$-supercritical case; otherwise, it is referred to as the mass mixed case.

In the purely $L^2$-supercritical case where $N=3, \alpha=\beta=2, p=q=4$, Bartsch, Jeanjean and Soave \cite[Theorem 1.1 and Theorem 1.2]{Bartsch2016} established the existence of positive normalized solutions for $0<\nu<\nu_1$ and $\nu>\nu_2$. Here, $\nu_1,\nu_2>0$ are implicitly defined by
$$\max\left\{\frac{1}{a^2\mu_1^2},\frac{1}{b^2\mu_2^2}\right\}=\frac{1}{a^2(\mu_1+\nu_1)^2}+\frac{1}{b^2(\mu_2+\nu_1)^2}$$
and
$$\frac{(a^2+b^2)^3}{(\mu_1a^4+\mu_2b^4+2\nu_2a^2b^2)^2}=\min\left\{\frac{1}{a^2\mu_1^2},\frac{1}{b^2\mu_2^2}\right\}.$$
Clearly, the bounds $\nu_1,\nu_2$ depend on the masses $a,b>0$, and they satisfy
\beq\lab{eq:property-nu}
\nu_1\to0,\ \nu_2\to\infty \quad\text{as $\frac{a}{b}\to0$ or $\frac{a}{b}\to\infty$}.
\eeq
In particular, no value of $\nu$ exists for which the results in \cite{Bartsch2016} guarantee a solution across all masses. This raises a natural question: Can $\nu_1,\nu_2$ be improved, and what is the optimal range of $(a,b,\nu)$ for the existence of positive normalized solutions? Bartsch, Jeanjean, and Soave posed this as an open problem in \cite[Remark 1.3-(i)]{Bartsch2016}. Subsequently, Bartsch, Zhong, and Zou \cite{Bartsch2021} addressed this problem by developing a novel approach, distinct from the variational framework on the constrained manifold, to explore almost the best range for the existence of positive normalized solutions. Their work, now referred to in the literature as the global branch approach, provides a resolution to the open problem posed by Bartsch-Jeanjean-Soave.

For the general purely $L^2$-supercritical case where $2\leq N\leq 4$ with $\alpha>1,\beta>1$, and parameters satisfying $2+\frac{4}{N}<p,q,\alpha+\beta<2^*$, Bartsch and Jeanjean \cite[Theorem 2.4]{Bartsch2018} established existence results for $0<\nu<\nu_1$ and $\nu>\nu_2$. Here $\nu_1,\nu_2$ depend explicitly on the masses $a,b$. Recently, Jeanjean, Zhang, and Zhong \cite{Jeanjean2024a} obtained an improved existence range for $(a,b,\nu)$ using a variational approach constrained to $L^2$-balls. We now detail this methodology. For notational convenience, we denote throughout this work $E=H^1(\R^N)\times H^1(\R^N)$ and $E_{rad}=H^1_{rad}(\R^N)\times H^1_{rad}(\R^N)$.

Within the variational approach constrained to $L^2$-balls, the following Pohozaev manifold plays a fundamental role in our analysis:
$$
\cp:=\{(u,v)\in E\setminus\{(0,0)\}:P(u,v)=0\},
$$
where
$$
P(u,v)=\|\nabla u\|_2^2+\|\nabla v\|_2^2-\frac{\mu_1}{p}\gamma_p\|u\|_p^p-\frac{\mu_2}{q}\gamma_q\|v\|_q^q-\gamma_{\alpha+\beta}\nu\int_{\R^N}|u|^{\alpha}|v|^{\beta}{\rm d}x
$$
with $\gamma_p=\frac{(p-2)N}{2}$. By the Pohozaev identity, any critical point of $I|_{T(a,b)}$ lies in $\cp\cap T(a,b)$. Consequently, demonstrating that $$\inf\limits_{\cp\cap T(a,b)}I(u,v)$$
is achieved at a critical point would naturally yield the existence of a ground state. To address this problem, we proceed in two steps:

{\bf Step 1:} Consider a relaxed problem.
For $a\geq0$, define the set
\beq\label{eq:240920-1}
D_a:=\{u\in H^1(\R^N):\|u\|_2^2\leq a\},
\eeq
and for any $a>0,b>0$, introduce
\beq\label{eq:240920-2}
\cp_{a,b}:=\cp\cap(D_a\times D_b).
\eeq
The relaxed problem is to seek a critical point $(u,v)\in E\setminus\{(0,0)\}$ of $I|_{D_a\times D_b}$ at the level
\beq\label{eq:240822-1}
C_{a,b}:=\inf\limits_{\cp_{a,b}}I(u,v).
\eeq
{\bf Step 2:} We demonstrate that any minimizer $(u,v)$ belongs to $T(a,b)$ and constitutes a critical point of $I|_{T(a,b)}$ at the level $\inf\limits_{\cp\cap T(a,b)}I(u,v)$.

Solutions of the form $(u,0)$ satisfying system \eqref{eq:240722-1} with $\|u\|_2^2=a$, or $(0,v)$ with $\|v\|_2^2=b$, are conventionally termed semi-trivial solutions. A key advantage of the variational framework on $L^2$-balls lies in enabling local analysis of $I$ restricted to $\cp_{a,b}$ near these semi-trivial solutions. This approach establishes the strict inequality
$$
C_{a,b}<\min\{I(u,0), I(0,v)\},
$$
see \cite[Lemma 7.3]{Jeanjean2024a}, which subsequently excludes semi-trivial solutions as weak limits of Palais-Smale sequences. In the present work, we adapt this methodology to establish existence of normalized solutions.

\subsection{Main results}
In this paper, we consider the mass mixed case that
\beq\label{condition:240901-2}
2<p<2+\frac{4}{N}<q<2^*,\alpha>1, \beta>1, r:=\alpha+\beta\leq2^*,
\eeq
where the Sobolev critical exponent is given by
$$
2^*=
\begin{cases}
+\infty, &N=1,2,\\
\frac{2N}{N-2},  &N\geq3.
\end{cases}
$$

For the Sobolev subcritical mass mixed case satisfying
\beq\label{condition:240901-1}
2<p<2+\frac{4}{N}<q<2^*,\alpha>1, \beta>1, r:=\alpha+\beta<2^*,
\eeq
Bartsch and Jeanjean \cite[Theorem 2.2]{Bartsch2018} established the existence of positive normalized solutions for certain ranges of $a,b$ when $\beta>2$; while Li and Zou \cite[Theorem 1.1]{Li2021a} proved existence of positive normalized ground states for system \eqref{eq:240722-1}-\eqref{eq:240722-2} under $\beta<2$. However, it is well known that for the mass mixed case, system typically exhibit solution multiplicity. Exploring this multiplicity constitutes a primary motivation for the present study.

Our first result establishes the existence of normalized ground state solutions for the Sobolev-subcritical mass-mixed case with $\beta=2$.
\bt\label{th:240829-1}
Assume $1\leq N\leq4$, \eqref{condition:240901-1} and $\beta=2$. For any $b>0$ and $\nu>\nu_{\mu_1,p,a,N,\alpha}$, there exists $a(b,\nu)>0$ such that for any $a<a(b,\nu)$, the system \eqref{eq:240722-1}-\eqref{eq:240722-2} admits a positive normalized ground state. Here $\nu_{\mu_1,p,a,N,\alpha}$ is defined by
\beq\label{eq:240829-2}
\nu_{\mu_1,p,a,N,\alpha}=\frac{1}{2}\mu_1^{\frac{\beta N-4}{N(p-2)-4}}\|U_p\|_2^{\frac{4(p-2-\alpha)}{N(p-2)-4}}a^{\frac{-2(p-2-\alpha)}{N(p-2)-4}}
\inf\limits_{h\in H^1(\R^N)\setminus\{0\}}\frac{\int_{\R^N}|\nabla h|^2{\rm d}x}{\int_{\R^N}U_p^{\alpha}h^2{\rm d}x},
\eeq
where $U_p$ denotes the unique positive solution to
$$
\begin{cases}
-\Delta u +u=u^{p-1} ~\hbox{in}~ \R^N,\\
u>0 ~\hbox{in}~ \R^N,\\
u(0)=\max\limits_{x\in\R^N}u(x) ~\hbox{and}~ u\in H^1(\R^N).
\end{cases}
$$
Furthermore, $\nu_{\mu_1,p,a,N,\alpha}=0$ when $N=1,2$.
\et
\br
\begin{itemize}
\item [(i)] In proving Theorem \ref{th:240829-1}, by Ekeland variational principle, one can obtain the existence of Palais-Smale sequence. As mentioned previously, to obtain the compactness, establishing compactness requires excluding semi-trivial solutions as weak limits of this sequence. Defining $m_p^{\mu_1}(a)$ and $m_q^{\mu_2}(b)$ as in Lemma \ref{le:240729-1}, the critical energy estimate
    \beq\label{eq:20241106-1}
    m(a,b)<\min\{m_p^{\mu_1}(a),m_q^{\mu_2}(b)\},
    \eeq
    (established in Lemma \ref{le:240830-1}) becomes essential. Crucially, the variational framework on $L^2$-balls enables existence of normalized ground states for system\eqref{eq:240722-1}-\eqref{eq:240722-2} for all $\beta\leq2$. Compared with \cite[Theorem 1.1]{Li2021a}, for the technical reasons, they only obtain the estimation \eqref{eq:20241106-1} for $\beta<2$ (see \cite[Lemma 3.4]{Li2021a}). So this also highlights a key advantage of the $L^2$-balls constrained approach.
\item [(ii)] The threshold $\nu_{\mu_1,p,a,N,\alpha}$ originates in \cite[(2.5)]{Jeanjean2024a}, the authors of \cite{Jeanjean2024a} show that \eqref{eq:240829-2} is equivalent to
    $$
    \nu_{\mu_1,p,a,N,\alpha}=\frac{1}{2}\inf\limits_{h\in H^1(\R^N)\setminus\{0\}}\frac{\int_{\R^N}|\nabla h|^2{\rm d}x}{\int_{\R^N}|z_{p,\mu_1,a}|^{\alpha}|h|^2{\rm d}x}.
    $$
    Here $z_{p,\mu_1,a}$ denotes the unique solution to \eqref{eq:240726-1} with parameters $\mu=\mu_1$ and $\eta=p$.
\end{itemize}
\er
In fact, the positive normalized ground state in Theorem \ref{th:240829-1} is also a local minimizer. Based on the existence of this local minimizer and the fact that $I$ is unbounded below on $T(a,b)$, it is natural to expect a second solution of mountain pass type solution. Building upon Theorem \ref{th:240829-1} and \cite[Theorem 1.1]{Li2021a}, our next theorem establishes the existence of multiple normalized solutions for all $\beta\leq2$.
\bt\label{th:240829-2}
Suppose $2\leq N\leq4$ and that \eqref{condition:240901-1} holds. Then the following statements hold:
\begin{itemize}
\item [(i)] If $\alpha,\beta<2$, then for any $\nu>0$ and $b>0$, there exists $a(b,\nu)>0$ such that for any $a<a(b,\nu)$, the system \eqref{eq:240722-1}-\eqref{eq:240722-2} admits a second solution in $E_{rad}$, which is of mountain pass type;
\item [(ii)] If $\alpha<2,\beta=2$, then for any $\nu>\nu_{\mu_1,p,a,N,\alpha}$ and $b>0$, there exists $a(b,\nu)>0$ such that for any $a<a(b,\nu)$, the conclusion in (i) holds;
\item [(iii)] If $\alpha=2,\beta<2$, then for any $\nu>\nu_{\mu_2,q,b,N,\beta}$ and $b>0$, there exists $a(b,\nu)>0$ such that for any $a<a(b,\nu)$, the conclusion in (i) holds;
\item [(iv)] If $\alpha=2,\beta=2$, then for any $\nu>\max\{\nu_{\mu_1,p,a,N,\alpha},\nu_{\mu_2,q,b,N,\beta}\}$ and $b>0$, there exists $a(b,\nu)>0$ such that for any $a<a(b,\nu)$, the conclusion in (i) holds.
\end{itemize}
In particular, $\nu_{\mu_1,p,a,N,\alpha}=0$ and $\nu_{\mu_2,q,b,N,\beta}=0$ when $N=2$.
\et
Nowadays, the research on Schr\"odinger system involving Sobolev critical exponent is also a hot topic. Beyond their significance in applications, not negligible reasons of the mathematicians interest for such problems are their stimulating and challenging difficulties coming from the lack of compactness, due to the presence of the limiting exponents for the Sobolev embedding theorems. For the scalar equation case, we can refer to \cite{Soave2020a,Jeanjean2022,Wei2022,Vicentiu2023a,Li2024,Li2021,Alves2022}. Significant progress has also been made for system \eqref{eq:240722-1}-\eqref{eq:240722-2}. Specifically, in the purely $L^2$-supercritical case where $p=q=2^*$ and $\alpha+\beta>\bar p$, Bartsch, Li and Zou \cite{Bartsch2023} established the existence of a normalized ground state solution provided $\nu>\nu_1$. Here
$$
\begin{cases}
\nu_1=0, ~\hbox{if}~ N=4 ~\hbox{or}~ N=3 ~\hbox{and}~ |\alpha-\beta|>2,\\
\nu_1>0, ~\hbox{if}~ N=3, \alpha\geq2,\beta\geq2.
\end{cases}
$$
In the mass mixed case where $p=q=2^*$, $\alpha+\beta<2+\frac{4}{N}$ and $N\in\{3,4\}$, Bartsch, Li and Zou \cite{Bartsch2023} also established the existence of a normalized ground state solution. They left the existence of a second normalized solution as an open problem (see \cite[Remark 1.3]{Bartsch2023}). Subsequently, Zhang, Zhong and Zhou provided a positive answer to this problem in \cite{Vicentiu2023}. In another mass mixed case with $2<p,q<\bar p$ and $\alpha+\beta=2^*$, He, Guo, Shuai and Zhong \cite{He2024} demonstrated that system \eqref{eq:240722-1}-\eqref{eq:240722-2} admits multiple positive normalized solutions for all $N\geq3$. In \cite{He2024}, the authors ingeniously constructed a special mountain pass geometry and precisely estimated the mountain pass level to restore the compactness of the Palais-Smale sequence.

Motivated by these results, in this paper, we focus on the following Sobolev critical case:
\beq\label{condition:240722-1}
2<p<2+\frac{4}{N}<q<2^*,\alpha>1, \beta>1, r=\alpha+\beta=2^*.
\eeq
We aim to establish the existence of multiple positive normalized solutions for system \eqref{eq:240722-1}-\eqref{eq:240722-2}.
\bt\label{th:240722-1}
Assume \eqref{condition:240722-1} holds and $N\in\{3,4\}$. Then the following statements hold:
\begin{itemize}
\item[(i)] If $1<\beta<2$, then for any $\nu>0$ and $b>0$, there exists $a(b,\nu)>0$ such that for any $a<a(b,\nu)$, $I$ admits a positive, radially symmetric, and decreasing local minimizer $(u,v)$ on $T(a,b)$. The corresponding Lagrange multipliers satisfy $\lambda_1,\lambda_2>0$. Moreover, the local minimizer $(u,v)$ is indeed a normalized ground state.
\item[(ii)] If $\beta=2$, then for any $\nu>\nu_{\mu_1,p,a,N,2^*-2}$ and $b>0$, there exists $a(b,\nu)>0$ such that for any $a<a(b,\nu)$, the conclusions in (i) also holds.
Here $\nu_{\mu_1,p,a,N,2^*-2}$ is defined by \eqref{eq:240829-2}.
\end{itemize}
\et

\bt\label{th:240722-2}
Assume the hypotheses of Theorem \ref{th:240722-1} holds. Let $\nu$ and $(u,v)$ be given in Theorem \ref{th:240722-1}. Then there exists $b_0>0$ such that for any $b<b_0$, we can find $a(b,\nu)>0$ such that for any $a<a(b,\nu)$, the system \eqref{eq:240722-1}-\eqref{eq:240722-2} admits a second solution $(u_0,v_0)$ in $E_{rad}$, which is of mountain pass type and satisfies
$$
I(u,v)<0<I(u_0,v_0)<I(u,v)+\frac{2}{N-2}\nu^{-\frac{N-2}{2}}\alpha^{-\frac{(N-2)\alpha}{4}}\beta^{-\frac{(N-2)\beta}{4}}S^{\frac{N}{2}}.
$$
Moreover, the corresponding Lagrange multipliers are positive: $\bar\lambda_1,\bar\lambda_2>0$.
\et
\br
Compared to Theorem \ref{th:240829-2}, we cannot simultaneously require $\alpha\leq2$ and $\beta\leq2$ since $\alpha+\beta=2^*$.
So, to exclude the semi-trivial solution for technical reasons, we restrict to sufficiently small $b>0$. By Lemma \ref{le:240727-2}, the energy of semi-trivial solution $m_q^{\mu_2}(b)$ is decreasing in $b$ with $\lim\limits_{b\to0^+}m_q^{\mu_2}(b)=+\infty$. Consequently, there exists $b_0>0$ such that for all $0<b<b_0$,
$$
m_q^{\mu_2}(b)\geq \frac{2}{N-2}\nu^{-\frac{N-2}{2}}\alpha^{-\frac{(N-2)\alpha}{4}}\beta^{-\frac{(N-2)\beta}{4}}S^{\frac{N}{2}}.
$$
Combining this with Corollary \ref{le:240801-2}, we obtain $0<M_{a,b}<m_q^{\mu_2}(b)$. This ensures that the weak limit of the Palais-Smale-Pohozaev (PSP) sequence cannot be a semi-trivial solution.
\er
The paper is structured as follows. Section \ref{preliminary} presents preliminary results. Section \ref{sec:240901-1} establishes essential properties of the Pohozaev manifold for subsequent analysis. Section \ref{sec:subcritical} addresses the Sobolev subcritical case, containing the proofs of Theorems \ref{th:240829-1} and \ref{th:240829-2}. Section \ref{sec:240829-1} focuses on the Sobolev critical case and contains the proofs of Theorems \ref{th:240722-1} and \ref{th:240722-2}.

\s{Preliminaries}\label{preliminary}
We begin by recalling some classical facts required for the subsequent analysis.
\bl\cite[Lemma 3.1]{Bartsch2018}\label{le:240729-1}
Assume that $\eta\in(2,2^*)\setminus\{2+\frac{4}{N}\}$ and let $\mu>0$ be given. For any $a>0$, there exists a unique couple $(\lambda_a,u_a)\in \R^+\times H^1_{rad}(\R^N)$ solving
\beq\label{eq:240726-1}
\begin{cases}
-\Delta u+\lambda u=\mu |u|^{\eta-2}u, u\in H^1(\R^N),\\
u>0,\\
\|u\|_2^2=a.
\end{cases}
\eeq
Moreover, $u_a$ corresponds to the least energy level $m_{\eta}^{\mu}(a)$ of the functional $I: H^1(\R^N)\to \R$ defined by
$$I(u)=\frac{1}{2}\|\nabla u\|_2^2-\frac{\mu}{\eta}\|u\|_{\eta}^{\eta}$$
constrained to the $L^2$-sphere $S_a$.
If $\eta\in(2,2+\frac{4}{N})$, then $m_{\eta}^{\mu}(a)<0$ for all $a>0$.
If $\eta\in(2+\frac{4}{N},2^*)$, then $m_{\eta}^{\mu}(a)>0$ for all $a>0$.
\el

\bl\cite[Lemma 3.2]{Jeanjean2024a}\label{le:240727-2}
$m_{\eta}^{\mu}(a)$ is continuous and decreases strictly respect to $a\in\R^+$. Furthermore,
\begin{itemize}
\item [(i)] If $2<\eta<2+\frac{4}{N}$, then
$$\lim\limits_{a\to 0^+}m_{\eta}^{\mu}(a)=0 ~\hbox{and}~ \lim\limits_{a\to +\infty}m_{\eta}^{\mu}(a)=-\infty.$$
\item [(ii)] If $2+\frac{4}{N}<\eta<2^*$, then
$$\lim\limits_{a\to 0^+}m_{\eta}^{\mu}(a)=+\infty ~\hbox{and}~ \lim\limits_{a\to +\infty}m_{\eta}^{\mu}(a)=0.$$
\end{itemize}
\el

As mentioned earlier, the Pohozaev manifold $\cp$ plays a crucial role in the proof. We retain the previous definitions for:
\begin{itemize}
    \item the $L^2$-ball $\mathcal{D}_a$,
    \item the Pohozaev manifold $\mathcal{P}_{a,b}$,
    \item the infimum $C_{a,b}$,
\end{itemize}
as given in \eqref{eq:240920-1}--\eqref{eq:240822-1}.

To analyze the minimization problem \eqref{eq:240822-1}, we introduce a dilation operator preserving the $L^2$-norm. For $u\in S_a$ and $t\in\R^+$, define
\beq\label{eq:240731-4}
t\star u(x)=t^{\frac{N}{2}}u(tx) ~\hbox{for a.e.}~ x\in\R^N.
\eeq
Consider the fiber map associated to $(u,v)$:
\beq\label{eq:240824-1}
\Phi_{(u,v)}(t):=I(t\star(u,v))=\frac{1}{2}(\|\nabla u\|_2^2+\|\nabla v\|_2^2)t^2-\frac{\mu_1}{p}\|u\|_p^pt^{\gamma_p}-\frac{\mu_2}{q}\|v\|_q^qt^{\gamma_q}-\nu\int_{\R^N}|u|^{\alpha}|v|^{\beta}{\rm d}xt^{\gamma_{\alpha+\beta}}.
\eeq
Direct computation shows that $t\Phi'_{(u,v)}(t)=P(t\star(u,v))$ and
$$
\cp_{a,b}=\{(u,v)\in D_a\times D_b\setminus\{(0,0)\}:\Phi'_{(u,v)}(1)=0\}.
$$
We decompose this manifold into the disjoint union $\cp_{a,b}=\cp_{a,b}^{+}\cup\cp_{a,b}^{0}\cup\cp_{a,b}^{-}$, where
\begin{align}\label{eq:240829-1}
\cp_{a,b}^{+}:=\{(u,v)\in\cp_{a,b}:\Phi''_{(u,v)}(1)>0\},\nonumber\\
\cp_{a,b}^{0}:=\{(u,v)\in\cp_{a,b}:\Phi''_{(u,v)}(1)=0\},\nonumber\\
\cp_{a,b}^{-}:=\{(u,v)\in\cp_{a,b}:\Phi''_{(u,v)}(1)<0\}.
\end{align}

To investigate fundamental properties of the Pohozaev manifold $\cp$, we first establish the following lemmas.

\noindent\textbf{Case $r := \alpha + \beta = 2^*$:}
Recall the Gagliardo-Nirenberg inequality for $N\geq3$ and $2<p\leq2^*$:
\beq\label{eq:240723-1}
\|u\|_p^p\leq C_{N,p}\|\nabla u\|_2^{\gamma_p}\|u\|_2^{p-\gamma_p},\quad \forall u\in H^1(\R^N).
\eeq
Here $\gamma_p=\frac{(p-2)N}{2}$. When $p=2^*$, this reduces to the Sobolev inequality:
\beq\label{eq:240723-2}
\|\nabla u\|_2^2\geq S\|u\|_{2^*}^{2}, \quad \forall u\in D_0^{1,2}(\R^N).
\eeq
Applying \eqref{eq:240723-1}, \eqref{eq:240723-2} and H\"{o}lder inequality, for any $(0,0)\neq(u,v)\in D_a\times D_b$, we derive
\begin{align}\label{eq:240724-3}
I(u,v)=&\frac{1}{2}(\|\nabla u\|_2^2+\|\nabla v\|_2^2)-\frac{\mu_1}{p}\|u\|_p^p-\frac{\mu_2}{q}\|v\|_q^q-\nu\int_{\R^N}|u|^{\alpha}|v|^{\beta}{\rm d}x\nonumber\\
\geq &\frac{1}{2}(\|\nabla u\|_2^2+\|\nabla v\|_2^2)-\frac{\mu_1}{p}C_{N,p}\|\nabla u\|_2^{\gamma_p}\|u\|_2^{p-\gamma_p}-\frac{\mu_2}{q}C_{N,q}\|\nabla v\|_2^{\gamma_q}\|v\|_2^{q-\gamma_q}-\nu\|u\|_{2^*}^{\alpha}\|v\|_{2^*}^{\beta}\nonumber\\
>&\frac{1}{2}(\|\nabla u\|_2^2+\|\nabla v\|_2^2)-\frac{\mu_1}{p}C_{N,p}a^{\frac{p-\gamma_p}{2}}(\|\nabla u\|_2^2+\|\nabla v\|_2^2)^{\frac{\gamma_p}{2}}\nonumber\\
&-\frac{\mu_2}{q}C_{N,q}b^{\frac{q-\gamma_q}{2}}(\|\nabla u\|_2^2+\|\nabla v\|_2^2)^{\frac{\gamma_q}{2}}-\nu S^{-\frac{2^*}{2}}(\|\nabla u\|_2^2+\|\nabla v\|_2^2)^{\frac{2^*}{2}}\nonumber\\
=:&h((\|\nabla u\|_2^2+\|\nabla v\|_2^2)^{\frac{1}{2}}),
\end{align}
where
\beq\label{eq:240723-4}
h(t)=\frac{1}{2}t^2-D_1t^{\gamma_p}-D_2t^{\gamma_q}-\nu S^{-\frac{2^*}{2}}t^{2^*}, \quad t\in (0,+\infty)
\eeq
with $D_1=\frac{\mu_1}{p}C_{N,p}a^{\frac{p-\gamma_p}{2}}$ and $D_2=\frac{\mu_2}{q}C_{N,q}b^{\frac{q-\gamma_q}{2}}$.

\bl\label{le:240724-1}
Assume condition \eqref{condition:240722-1} holds and $N\geq3$. Define the threshold parameter
\beq\label{eq:240723-3}
T_{a,b}:=\mu_2 b^{\frac{q-\gamma_q}{2}}(\mu_1 a^{\frac{p-\gamma_p}{2}})^{\frac{\gamma_q-2}{2-\gamma_p}}+\nu(\mu_1 a^{\frac{p-\gamma_p}{2}})^{\frac{2^*-2}{2-\gamma_p}}.
\eeq
There exists $\alpha_1>0$ such that if $T_{a,b}<\alpha_1$, then $h(t)$ admits exactly two critical points: a local minimum with the negative level, the other is a global maximum with the positive level. Furthermore, there exist unique $0<R_0<R_1$ such that $h(R_0)=h(R_1)=0$ and
$$
h(t) > 0 \quad \text{if and only if} \quad t \in (R_0, R_1).
$$
\el
\bp
By \eqref{eq:240723-4}, we have
$$
\begin{aligned}
h'(t)=&t-D_1\gamma_pt^{\gamma_p-1}-D_2\gamma_qt^{\gamma_q-1}-2^*\nu S^{-\frac{2^*}{2}}t^{2^*-1}\\
=&t^{\gamma_p-1}(t^{2-\gamma_p}-D_1\gamma_p-D_2\gamma_qt^{\gamma_q-\gamma_p}-2^*\nu S^{-\frac{2^*}{2}}t^{2^*-\gamma_p})\\
=&:t^{\gamma_p-1}(g(t)-D_1\gamma_p),
\end{aligned}
$$
where
$$g(t)=t^{2-\gamma_p}-D_2\gamma_qt^{\gamma_q-\gamma_p}-2^*\nu S^{-\frac{2^*}{2}}t^{2^*-\gamma_p}.$$
Direct computation yields the derivative of $g$:
$$
\begin{aligned}
g'(t)=&(2-\gamma_p)t^{1-\gamma_p}-D_2\gamma_q(\gamma_q-\gamma_p)t^{\gamma_q-\gamma_p-1}-2^*\nu S^{-\frac{2^*}{2}}({2^*-\gamma_p})t^{2^*-\gamma_p-1}\\
=&t^{1-\gamma_p}(2-\gamma_p-D_2\gamma_q(\gamma_q-\gamma_p)t^{\gamma_q-2}-2^*\nu S^{-\frac{2^*}{2}}({2^*-\gamma_p})t^{2^*-2}).
\end{aligned}
$$
Setting $g'(t)=0$ gives an equation whose unique positive solution $\bar t>0$ satisfies
\beq\label{eq:240724-1}
(2-\gamma_p){\bar t}^2=D_2\gamma_q(\gamma_q-\gamma_p){\bar t}^{\gamma_q}+2^*\nu S^{-\frac{2^*}{2}}({2^*-\gamma_p}){\bar t}^{2^*}.
\eeq
Observe that if $g(\bar t)>D_1\gamma_p$ and additionally $h(\bar t)>0$, then $h'(t)=0$ has exactly two solutions. In this case, there exist $0<R_0<R_1$ such that $h(R_0)=h(R_1)=0$ and $h(t)>0$ if and only if $t\in (R_0,R_1)$. Therefore, it suffices to establish the conditions $g(\bar t)>D_1\gamma_p$ and $h(\bar t)>0$, i.e.,
\beq\label{eq:240724-2}
\begin{cases}
{\bar t}^2>D_1\gamma_p {\bar t}^{\gamma_p}+D_2\gamma_q {\bar t}^{\gamma_q}+2^*\nu S^{-\frac{2^*}{2}}{\bar t}^{2^*},\\
\frac{1}{2}{\bar t}^2>D_1{\bar t}^{\gamma_p}+D_2{\bar t}^{\gamma_q}+\nu S^{-\frac{2^*}{2}}{\bar t}^{2^*}.
\end{cases}
\eeq
From \eqref{eq:240724-1}, we obtain the estimates:
\begin{align*}
&D_2\gamma_q {\bar t}^{\gamma_q}+2^*\nu S^{-\frac{2^*}{2}}{\bar t}^{2^*}<\frac{2-\gamma_p}{\gamma_q-\gamma_p}{\bar t}^2\\
&D_2{\bar t}^{\gamma_q}+\nu S^{-\frac{2^*}{2}}{\bar t}^{2^*}<\frac{1}{2}(D_2\gamma_q {\bar t}^{\gamma_q}+2^*\nu S^{-\frac{2^*}{2}}{\bar t}^{2^*})<\frac{2-\gamma_p}{2(\gamma_q-\gamma_p)}{\bar t}^2.\\
\end{align*}
This implies the following inequalities:
$$
\begin{cases}
D_1\gamma_p {\bar t}^{\gamma_p}+D_2\gamma_q {\bar t}^{\gamma_q}+2^*\nu S^{-\frac{2^*}{2}}{\bar t}^{2^*}-{\bar t}^2<\left(\frac{2-\gamma_q}{\gamma_q-\gamma_p}+D_1\gamma_p{\bar t}^{\gamma_p-2}\right){\bar t}^2,\\
D_1{\bar t}^{\gamma_p}+D_2{\bar t}^{\gamma_q}+\nu S^{-\frac{2^*}{2}}{\bar t}^{2^*}-\frac{1}{2}{\bar t}^2<\left(\frac{1}{2}\frac{2-\gamma_q}{\gamma_q-\gamma_p}+D_1{\bar t}^{\gamma_p-2}\right){\bar t}^2.
\end{cases}
$$
Define $\bar s=\left(2D_1\frac{\gamma_q-\gamma_p}{\gamma_q-2}\right)^{\frac{1}{2-\gamma_p}}$. If $\bar t>\bar s$, then \eqref{eq:240724-2} holds. To ensure $\bar t>\bar s$, by \eqref{eq:240724-1} it suffices to verify
\beq\label{eq:crit_cond}
(2-\gamma_p){\bar s}^2>D_2\gamma_q(\gamma_q-\gamma_p){\bar s}^{\gamma_q}+2^*\nu S^{-\frac{2^*}{2}}({2^*-\gamma_p}){\bar s}^{2^*}.
\eeq
Substituting $\bar s=\left(2D_1\frac{\gamma_q-\gamma_p}{\gamma_q-2}\right)^{\frac{1}{2-\gamma_p}}$ into \eqref{eq:crit_cond} yields
\beq\label{eq:substituted}
(2-\gamma_p)>D_2\gamma_q(\gamma_q-\gamma_p)\left(2D_1\frac{\gamma_q-\gamma_p}{\gamma_q-2}\right)^{\frac{\gamma_q-2}{2-\gamma_p}}+2^*\nu S^{-\frac{2^*}{2}}({2^*-\gamma_p})\left(2D_1\frac{\gamma_q-\gamma_p}{\gamma_q-2}\right)^{\frac{2^*-2}{2-\gamma_p}}.
\eeq
Now substitute $D_1=\frac{\mu_1}{p}C_{N,p}a^{\frac{p-\gamma_p}{2}}$ and $D_2=\frac{\mu_2}{q}C_{N,q}b^{\frac{q-\gamma_q}{2}}$ into \eqref{eq:substituted}:
\begin{align}\label{eq:20250604-1}
2-\gamma_p>
&\frac{\mu_2C_{N,q}}{q}b^{\frac{q-\gamma_q}{2}}\gamma_q(\gamma_q-\gamma_p)\left(\frac{2\mu_1C_{N,p}}{p}a^{\frac{p-\gamma_p}{2}}\frac{\gamma_q-\gamma_p}{\gamma_q-2}\right)^{\frac{\gamma_q-2}{2-\gamma_p}} \nonumber\\
+&2^*\nu S^{-\frac{2^*}{2}}(2^*-\gamma_p)\left(\frac{2\mu_1C_{N,p}}{p}a^{\frac{p-\gamma_p}{2}}\frac{\gamma_q-\gamma_p}{\gamma_q-2}\right)^{\frac{2^*-2}{2-\gamma_p}}\nonumber\\
>&\min\left\{\frac{C_{N,q}\gamma_q(\gamma_q-\gamma_p)}{q}\left(\frac{2C_{N,p}(\gamma_q-\gamma_p)}{p(\gamma_q-2)}\right)^{\frac{\gamma_q-2}{2-\gamma_p}},
\right.\nonumber\\
&\left. 2^*S^{-\frac{2^*}{2}}(2^*-\gamma_p)\left(\frac{2C_{N,p}(\gamma_q-\gamma_p)}{p(\gamma_q-2)}\right)^{\frac{2^*-2}{2-\gamma_p}}
\right\}T_{a,b}.
\end{align}
Define $$\alpha_1:=\frac{2-\gamma_p}{\min\left\{\frac{C_{N,q}\gamma_q(\gamma_q-\gamma_p)}{q}\left(\frac{2C_{N,p}(\gamma_q-\gamma_p)}{p(\gamma_q-2)}\right)^{\frac{\gamma_q-2}{2-\gamma_p}},
2^*S^{-\frac{2^*}{2}}(2^*-\gamma_p)\left(\frac{2C_{N,p}(\gamma_q-\gamma_p)}{p(\gamma_q-2)}\right)^{\frac{2^*-2}{2-\gamma_p}}\right\}}.$$
Then \eqref{eq:20250604-1} is equivalent to $T_{a,b}<\alpha_1$, which completes the proof.
\ep
\bl
Under the assumptions of Lemma \ref{le:240724-1}, for any $a,b>0$, the following holds:
\beq\label{eq:240825-1}
R_0<R_1<\left(\frac{1}{2\nu}S^{\frac{2^*}{2}}\right)^{\frac{1}{2^*-2}}.
\eeq
Furthermore, denoting $R_0=R_0(a,b)$ and $R_1=R_1(a,b)$, we have
\begin{align*}
R_0\to0 ~\hbox{and}~ R_1\to \left(\frac{1}{2\nu}S^{\frac{2^*}{2}}\right)^{\frac{1}{2^*-2}} ~\hbox{as}~ a,b\to0.
\end{align*}
\el
\bp
Recall from \eqref{eq:240723-4} that we denote the function $h(t)$ by $h_{a,b}(t)$. Observe that $h_{a,b}$ is continuous and monotonically decreasing with respect to $a,b\geq0$. When $a=b=0$, we have
$$
h_0(t)=\frac{1}{2}t^2-\nu S^{-\frac{2^*}{2}}t^{2^*}=t^2\left(\frac{1}{2}-\nu S^{-\frac{2^*}{2}}t^{2^*-2}\right)
$$
and
$$
h'_0(t)=t\left(1-2^*\nu S^{-\frac{2^*}{2}}t^{2^*-2}\right).
$$
Direct computation shows that $h_0(t)$ increases in $\left(0,\left(\frac{1}{2^*\nu}S^{\frac{2^*}{2}}\right)^{\frac{1}{2^*-2}}\right)$, decreases in $\left(\left(\frac{1}{2^*\nu}S^{\frac{2^*}{2}}\right)^{\frac{1}{2^*-2}},+\infty\right)$ and has a unique zero point at $t=\left(\frac{1}{2\nu}S^{\frac{2^*}{2}}\right)^{\frac{1}{2^*-2}}$.
For any $a>0,b>0$, the monotonicity of $h_{a,b}(t)$ with respect to $a$ and $b$ implies $h_0(R_0)>h_{a,b}(R_0)=0$ and $h_0(R_1)>h_{a,b}(R_1)=0$, which establishes \eqref{eq:240825-1}.

Since $R_0(a,b)$ and $R_1(a,b)$ are bounded for $a,b>0$, we may assume that $R_0(a,b)\to\eta_1$ and $R_1(a,b)\to\eta_2$ as $a,b\to0$. By Lemma \ref{le:240724-1}, $h(R_0)=0$ and $h'(R_0)>0$. Taking $a\to0,b\to0$, we obtain
\begin{align*}
h_0(\eta_1)=0 ~\hbox{and}~ h'_0(\eta_1)\geq0.
\end{align*}
So $\eta_1=0$. Similarly, from $h(R_1)=0$ and $h'(R_1)<0$, we obtain
\begin{align*}
h_0(\eta_2)=0 ~\hbox{and}~ h'_0(\eta_2)\leq0,
\end{align*}
which implies $\eta_2=\left(\frac{1}{2\nu}S^{\frac{2^*}{2}}\right)^{\frac{1}{2^*-2}}$.
\ep

\bl\label{le:240921-1}
Let $R_0=R_0(a,b)$ and $R_1=R_1(a,b)$ be defined as in Lemma \ref{le:240724-1}. Then $R_0(a,b)$ is increasing in each variable $a$ and $b$ for $a, b\geq0$, while $R_1(a,b)$ is decreasing in each variable $a$ and $b$ for $a, b\geq0$.
\el
\bp
Since $h_{a,b}(t)$ is jointly decreasing in $a$ and $b$ for $a, b\geq0$, it follows that for any $0\leq a_1\leq a$, $0\leq b_1\leq b$,
$$
h_{a_1,b_1}(R_0)\geq h_{a,b}(R_0)=0 ~\hbox{and}~ h_{a_1,b_1}(R_1)\geq h_{a,b}(R_1)=0.$$
Consequently, $R_0(a_1,b_1)\leq R_0(a,b)$ and $R_1(a_1,b_1)\geq R_1(a,b)$.
\ep

\bl\label{le:240824-1}
Let $2<p<2+\frac{4}{N}<q< 2^*$ and $l(t)=at^2-bt^{2^*}-ct^{\gamma_p}-dt^{\gamma_q}$, here $a, b, c, d>0, t\in(0,+\infty)$. Then $l(t)$ has at most two critical points.
\el
\bp
This result is obtained by a direct computation and we omit the details.
\ep

{\bf Case $r := \alpha + \beta < 2^*$:}
This case can be handled by referring to \cite[Section 3]{Li2021a}. To avoid redundancy, we omit the computational details and state the relevant results directly..

Adopting the notations from the case $\alpha+\beta=2^*$ and following \cite[(2.3)-(2.6)]{Li2021a}, we have, for any $(0,0)\neq(u,v)\in D_a\times D_b$,
\beq\label{eq:240901-1}
\begin{aligned}
I(u,v)>h((\|\nabla u\|_2^2+\|\nabla v\|_2^2)^{\frac{1}{2}}),
\end{aligned}
\eeq
where
$$
h(t)=\frac{1}{2}t^2-D_1t^{\gamma_p}-D_2t^{\gamma_q}-\nu D_3t^{\gamma_r}, t\in (0,+\infty)
$$
with the constants defined by
$$D_1=\frac{\mu_1}{p}C_{N,p}a^{\frac{p-\gamma_p}{2}}, D_2=\frac{\mu_2}{q}C_{N,q}b^{\frac{q-\gamma_q}{2}} ~\hbox{and}~ D_3=(\frac{\max\{\alpha,\beta\}}{r})^{\frac{\gamma_r}{2}}C_{N,r}a^{\frac{\alpha(r-\gamma_r)}{2r}}b^{\frac{\beta(r-\gamma_r)}{2r}}.$$
Furthermore, we note that Lemma \ref{le:240724-1} remains valid in this setting; for details we refer to \cite[Lemma 3.1]{Li2021a}.
\bl\cite[ Remark 3.1]{Li2021a}\label{le:240901-1}
Let $l(t)=at^2-bt^{\gamma_r}-ct^{\gamma_p}-dt^{\gamma_q}$ with $2<p<2+\frac{4}{N}<q\leq2^*$, $2<r<2^*$ and $a,b,c,d>0$. Then, $l(t)$ has at most two critical points in $(0,\infty)$.
\el

\s{Some properties of Pohozaev manifold $\cp$}\label{sec:240901-1}

In this section, for ease of reference, we recall that
$$
\widetilde T_{a, b}=
\begin{cases}
a^{\alpha\left(1-\frac{\gamma_r}{r}\right)} b^{\beta\left(1-\frac{\gamma_r}{r}\right)} \nu\left(\mu_2 b^{q-\gamma_q}\right)^{\frac{2-\gamma_r}{\gamma_q-2}}+\mu_1 a^{p-\gamma_p}\left(\mu_2 b^{q-\gamma_q}\right)^{\frac{2-\gamma_p}{\gamma_q-2}}, & \text { if } r<2+\frac{4}{N}, \\
\min \left\{a^{\alpha\left(1-\frac{\gamma_r}{r}\right)} b^{\beta\left(1-\frac{\gamma_r}{r}\right)} \nu,\left(\mu_1 a^{p-\gamma_p}\right)^{\frac{1}{2-\gamma_p}}\left(\mu_2 b^{q-\gamma_q}\right)^{\frac{1}{\gamma_q-2}}\right\}, & \text { if } r=2+\frac{4}{N}, \\
a^{\alpha\left(1-\frac{\gamma_r}{r}\right)} b^{\beta\left(1-\frac{\gamma_r}{r}\right)} \nu\left(\mu_1 a^{p-\gamma_p}\right)^{\frac{\gamma_r-2}{2-\gamma_p}}+\mu_2 b^{q-\gamma_q}\left(\mu_1 a^{p-\gamma_p}\right)^{\frac{\gamma_q-2}{2-\gamma_p}}, & \text { if } r>2+\frac{4}{N}.
\end{cases}
$$
and
$$
T_{a,b}:=\mu_2 b^{\frac{q-\gamma_q}{2}}(\mu_1 a^{\frac{p-\gamma_p}{2}})^{\frac{\gamma_q-2}{2-\gamma_p}}+\nu(\mu_1 a^{\frac{p-\gamma_p}{2}})^{\frac{2^*-2}{2-\gamma_p}}.
$$

\bl\label{le:240829-1}
Suppose $N\geq3$ and condition \eqref{condition:240722-1} holds, or $N\geq1$ and condition \eqref{condition:240901-1} holds. Then there exist positive constants $\alpha_2>0$ and $\tilde\alpha_2>0$ such that $\cp_{a,b}^{0}=\emptyset$ whenever $T_{a,b}<\alpha_2$ or $\widetilde T_{a,b}<\tilde\alpha_2$, respectively.
\el
\bp
For the case $\alpha+\beta<2^*$, the result follows directly from \cite[Lemma 3.2]{Li2021a}. Thus we focus on the critical case $\alpha+\beta=2^*$ and proceed by contradiction.

Suppose to the contrary that for every $\alpha_2>0$ with $T_{a,b}<\alpha_2$, there exists $(u,v)\in \cp_{a,b}$ satisfying $\Phi''_{(u,v)}(1)=0$. For such $(u,v)$, define the auxiliary function:
\begin{align*}
W(t):=(t-1)\Phi'_{(u,v)}(1)-\Phi''_{(u,v)}(1) \quad t\in \R.
\end{align*}
By assumption, $W(t)\equiv0$, which implies
$$
\begin{aligned}
0=&(t-2)(\|\nabla u\|_2^2+\|\nabla v\|_2^2)-(t-\gamma_p)\frac{\mu_1}{p}\gamma_p\|u\|_p^p\\
&-(t-\gamma_q)\frac{\mu_2}{q}\gamma_q\|v\|_q^q-(t-2^*)2^*\nu\int_{\R^N}|u|^{\alpha}|v|^{\beta}{\rm d}x.
\end{aligned}
$$
Set $\rho:=(\|\nabla u\|_2^2+\|\nabla v\|_2^2)^{\frac{1}{2}}$. Evaluating at $t=\gamma_q$ and applying \eqref{eq:240723-1} yields
$$
(\gamma_q-2)\rho^2\leq(\gamma_q-\gamma_p)\frac{\mu_1}{p}\gamma_p\|u\|_p^p\leq(\gamma_q-\gamma_p)\gamma_pD_1\rho^{\gamma_p},
$$
which implies
\beq\label{eq:240823-1}
\rho\leq\left(\frac{\gamma_q-\gamma_p}{\gamma_q-2}D_1\gamma_p\right)^{\frac{1}{2-\gamma_p}}.
\eeq

Now evaluate at $t=\gamma_p$. Using \eqref{eq:240723-1}, \eqref{eq:240823-1} and the Sobolev inequality, we obtain
$$
\begin{aligned}
(2-\gamma_p)=&(\gamma_q-\gamma_p)\frac{\mu_2}{q}\gamma_q\|v\|_q^q\rho^{-2}+(2^*-\gamma_p)2^*\nu\int_{\R^N}|u|^{\alpha}|v|^{\beta}{\rm d}x\rho^{-2}\\
\leq&(\gamma_q-\gamma_p)\gamma_qD_2\rho^{\gamma_q-2}+(2^*-\gamma_p)2^*\nu S^{-\frac{2^*}{2}}\rho^{2^*-2}\\
\leq&(\gamma_q-\gamma_p)\gamma_qD_2 \left(\frac{\gamma_q-\gamma_p}{\gamma_q-2}D_1\gamma_p\right)^{\frac{\gamma_q-2}{2-\gamma_p}}+
(2^*-\gamma_p)2^*\nu S^{-\frac{2^*}{2}} \left(\frac{\gamma_q-\gamma_p}{\gamma_q-2}D_1\gamma_p\right)^{\frac{2^*-2}{2-\gamma_p}}\\
\leq& C(p,q,N)(\mu_2 b^{\frac{q-\gamma_q}{2}}(\mu_1 a^{\frac{p-\gamma_p}{2}})^{\frac{\gamma_q-2}{2-\gamma_p}}+\nu(\mu_1 a^{\frac{p-\gamma_p}{2}})^{\frac{2^*-2}{2-\gamma_p}}),
\end{aligned}
$$
where $C(p,q,N)$ is a constant depending only on $p, q$ and $N$. This implies $T_{a,b}\geq {C(p,q,N)}^{-1}(2-\gamma_p)=:\alpha_2>0$, contradicting the assumption that $T_{a,b}$ can be arbitrarily small.
\ep

\bl\label{le:240824-2}
Suppose $N\geq3$ and condition \eqref{condition:240722-1} holds, or $N\geq1$ and condition \eqref{condition:240901-1} holds. Assume $T_{a,b}<\min\{\alpha_1,\alpha_2\}$ holds in the first case, or $\widetilde T_{a,b}<\min\{\tilde\alpha_1,\tilde\alpha_2\}$ in the second. Then for any $0<a_1\leq a$ and $0<b_1\leq b$, and for every $(u,v)\in T(a_1,b_1)\subset D_a\times D_b$, the function $\Phi_{(u,v)}(t)$ has exactly two critical points $0<s_{(u,v)}<t_{(u,v)}$. Moreover,
\begin{itemize}
\item [(i)] The scaling $s\star(u,v)$ belongs to $\cp_{a,b}^{+}$ if and only if $s=s_{(u,v)}$, and to $\cp_{a,b}^{-}$ if and only if $s=t_{(u,v)}$.
\item [(ii)] The first critical value satisfies:
\begin{align*}
\Phi_{(u,v)}(s_{(u,v)})=&\inf\left\{\Phi_{(u,v)}(t):t\in\left(0,\frac{R_0}{\left(\|\nabla u\|_2^2+\|\nabla v\|_2^2\right)^{\frac{1}{2}}}\right)\right\}\\
=&\inf\left\{\Phi_{(u,v)}(t):t\in\left(0,\frac{R_1}{\left(\|\nabla u\|_2^2+\|\nabla v\|_2^2\right)^{\frac{1}{2}}}\right)\right\}.
\end{align*}
\item [(iii)] The mappings $(u,v)\mapsto s_{(u,v)}$ and $(u,v)\mapsto t_{(u,v)}$ are of class $C^1$.
\end{itemize}
\el
\bp
We provide a unified proof covering both cases $\alpha+\beta<2^*$ and the case $\alpha+\beta=2^*$.

First, recall from \eqref{eq:240824-1} that $\Phi_{(u,v)}(0)=0$, while $\Phi_{(u,v)}(t)<0$ for sufficiently small $t>0$ and $\Phi_{(u,v)}(t)\to-\infty$ as $t\to+\infty$. Moreover, by \eqref{eq:240724-3} and Lemma \ref{le:240724-1}, for any $t\in\left(\frac{R_0}{\left(\|\nabla u\|_2^2+\|\nabla v\|_2^2\right)^{\frac{1}{2}}},\frac{R_1}{\left(\|\nabla u\|_2^2+\|\nabla v\|_2^2\right)^{\frac{1}{2}}}\right)$, we have $\Phi_{(u,v)}(t)>0$.

Second, Lemmas \ref{le:240824-1} and \ref{le:240901-1} imply that $\Phi_{(u,v)}(t)$ has at most two critical points in $(0,\infty)$. Combining this with the intermediate value theorem and the asymptotic behavior established above, we conclude that $\Phi_{(u,v)}(t)$ has exactly two critical points $s_{(u,v)},t_{(u,v)}$ satisfying
$$0<s_{(u,v)}<\frac{R_0}{\left(\|\nabla u\|_2^2+\|\nabla v\|_2^2\right)^{\frac{1}{2}}}<\frac{R_1}{\left(\|\nabla u\|_2^2+\|\nabla v\|_2^2\right)^{\frac{1}{2}}}.$$
This immediately yields conclusions (i) and (ii).

For (iii), consider the $C^1$ functional $\Psi(t,u,v):=\Phi'_{(u,v)}(t)$. At the critical points we have:
$$
\begin{aligned}
&\Psi(s_{(u,v)},u,v)=\Psi(t_{(u,v)},u,v)=0\\
&\partial_t \Psi(s_{(u,v)},u,v)=\Phi''_{(u,v)}(s_{(u,v)})>0\\
&\partial_t \Psi(t_{(u,v)},u,v)=\Phi''_{(u,v)}(t_{(u,v)})<0,
\end{aligned}
$$
Since $\partial_t \Psi \neq 0$ at both critical points, the Implicit Function Theorem guarantees that the mappings $(u,v) \mapsto s_{(u,v)}$ and $(u,v) \mapsto t_{(u,v)}$ are of class $C^1$.
\ep

\br\label{re:240918-1}
\begin{itemize}
\item [(i)] When $a=0$, we have $\cp_{0,b}=\{v\in D_b\setminus\{0\}: \Phi'_{(0,v)}(1)=0\}$; similarly, when $b=0$, $\cp_{a,0}=\{u\in D_a\setminus\{0\}: \Phi'_{(u,0)}(1)=0\}$. As established in the literature, $\cp_{0,b}=\cp_{0,b}^{-}$ and $\cp_{a,0}=\cp_{a,0}^{+}$.
\item [(ii)] For any non-trivial pair $(0,v)\in D_a\times D_b\setminus \{(0,0)\}$, the function $\Phi_{(0,v)}(t)$ has a unique critical point $t_{(0,v)}>0$ satisfying $t_{(0,v)}\star (0,v)\in\cp_{0,b}^{-}$. Analogously, for any non-trivial pair $(u,0)\in D_a\times D_b \setminus \{(0,0)\}$, the function $\Phi_{(u,0)}(t)$ has a unique critical point $s_{(u,0)}>0$ satisfying $s_{(u,0)}\star (u,0)\in\cp_{a,0}^{+}$.
\item [(iii)] Combining these observations with Lemma \ref{le:240824-2}, we conclude that for any $a>0$ and $b\geq0$, $C_{a,b}=\inf\limits_{\cp_{a,b}}I(u,v)=\inf\limits_{\cp_{a,b}^{+}}I(u,v)<0$.
\end{itemize}
\er

We now introduce, for each $R>0$, the set
$$
V_R:=\{(u,v)\in E:\|u\|_2^2\leq a,\|v\|_2^2\leq b,(\|\nabla u\|_2^2+\|\nabla v\|_2^2)^{\frac{1}{2}}\leq R\}.
$$
Setting $R=\frac{R_0+R_1}{2}$, we define
\beq\label{def:240901-1}
m(a,b):=\inf\limits_{V_R}I(u,v).
\eeq
The following property holds for this minimization problem:

\bl\label{le:240825-3}
For any $a>0$ and $b\geq0$, it holds that $C_{a,b}=m(a,b)$.
\el
\bp
First, for any $(u,v)\in\cp_{a,b}^{+}$, the conditions $\Phi'_{(u,v)}(1)=0$ and $I(u,v)=\Phi_{(u,v)}(1)<0$ imply that $(\|\nabla u\|_2^2+\|\nabla v\|_2^2)^{\frac{1}{2}}<R_0<R$. Consequently, $\cp_{(a,b)}^{+}\subset V_R$, and thus $C_{a,b}\geq m(a,b)$.

Second, for any $(u,v)\in V_{R}$, we have $1\in(0,\frac{R}{(\|\nabla u\|_2^2+\|\nabla v\|_2^2)^{\frac{1}{2}}}]$. By Lemma \ref{le:240824-2}-(ii), it follows that
$$
C_{a,b}\leq I(s_{(u,v)}\star(u,v))\leq I(u,v).
$$
Therefore, $C_{a,b}\leq m(a,b)$, which completes the proof.
\ep

\br\label{br:240729-1}
We remark that for any $a\geq0$ and $b\geq0$, the above $R=R(a,b)=\frac{R_0+R_1}{2}$ satisfies $R>R_0$. Moreover, we have the equality
$$m(a,b):=\inf\limits_{V_R}I(u,v)=\inf\limits_{V_{R_0}}I(u,v).$$
To see this, first observe that the function $h_{a,b}(t)$ is decreasing with respect to both $a\geq0$ and $b\geq0$. Now consider any pair $(u,v)$ in the set
$$\{(u,v)\in T(a_1,b_1): R_0\leq(\|\nabla u\|_2^2+\|\nabla v\|_2^2)^{\frac{1}{2}}\leq R, 0\leq a_1\leq a, 0\leq b_1\leq b,(a_1,b_1)\neq(0,0)\}.$$
For such $(u,v)$, we obtain
$$I(u,v)>h_{a_1,b_1}((\|\nabla u\|_2^2+\|\nabla v\|_2^2)^{\frac{1}{2}})\geq h_{a,b}((\|\nabla u\|_2^2+\|\nabla v\|_2^2)^{\frac{1}{2}})>0.$$
\er

\bl\label{le:240921-2}
The function $m(a,b)$ is decreasing in $a>0$ and $b\geq0$. That is,
$$m(a,b)\leq m(a_1,b_1) \quad \text{for all} \quad 0< a_1\leq a, 0\leq b_1\leq b.$$
\el
\bp
By Remark \ref{br:240729-1}, $m(a,b) = \inf\limits_{V_{R_0}} I(u,v)$. Lemma \ref{le:240921-1} states that $R_0(a,b)$ is increasing in $a \geq 0$ and $b \geq 0$. Therefore, for any $0< a_1 \leq a$ and $0 \leq b_1 \leq b$, we have $R_{0}(a_1,b_1) \leq R_0(a,b)$. Consequently, $V_{R_0(a_1,b_1)} \subset V_{R_0(a,b)}$, which implies $m(a_1,b_1) \geq m(a,b)$.
\ep

\s{The Sobolev subcritical case}\label{sec:subcritical}
As established in Section \ref{sec:240901-1}, when $\alpha + \beta < 2^*$, there exists a constant
$$c_0 = c_0(p, q, \alpha, \beta, N) = \min\{\tilde{\alpha}_1, \tilde{\alpha}_2\} > 0$$
ensuring that Lemmas \ref{le:240829-1}--\ref{le:240825-3} and Remark \ref{br:240729-1} hold whenever $\widetilde{T}_{a,b} < c_0$.
Recalling definition \eqref{eq:240829-1}, we define
\beq\label{def:240830-1}
Z_{a,b} := \inf\limits_{\mathcal{P}_{a,b}^{-}} I(u,v).
\eeq

\subsection{Estimation of $C_{a,b}$ and $Z_{a,b}$}
\begin{lemma}\label{le:240830-1}
Assume condition \eqref{condition:240901-1}, $N \geq 1$, and $\widetilde{T}_{a,b} < c_0$. Let $m_{p}^{\mu_1}(a)$ and $m_{q}^{\mu_2}(b)$ be defined as in Lemma \ref{le:240729-1}, with $\nu_{\mu_1,p,a,N,\alpha}$ and $\nu_{\mu_2,q,b,N,\beta}$ given by \eqref{eq:240829-2}. Then:
\begin{itemize}
\item[(i)] For $1 < \beta < 2$, $C_{a,b} < m_{p}^{\mu_1}(a)$. When $\beta = 2$, $C_{a,b} < m_{p}^{\mu_1}(a)$ holds whenever $\nu > \nu_{\mu_1,p,a,N,\alpha}$.
\item[(ii)] For $1 < \alpha < 2$, $Z_{a,b} < m_q^{\mu_2}(b)$. When $\alpha = 2$, $Z_{a,b} < m_q^{\mu_2}(b)$ holds whenever $\nu > \nu_{\mu_2,q,b,N,\beta}$.
\end{itemize}
\end{lemma}
\bp
The proof of (i) is completely similar to that of Lemma \ref{le:240825-1} and \cite[Lemma 7.3]{Jeanjean2024a}. We briefly prove (ii) below. For more details, see Lemma \ref{le:240825-1}.

Let $(v,\lambda)$ be the unique solution to equation \eqref{eq:240726-1} with $\mu=\mu_2$, $\eta=q$ and $\|v\|_2^2=b$. For any $h\in H^1(\R^N)$ with $\|h\|_2^2=1$, we have $(sh,v)\in D_a\times D_b$ provided $|s|\leq \sqrt{a}$. By Lemma \ref{le:240824-2}-(i) and Remark \ref{re:240918-1},
for each $s\geq0$, there exists a unique $t=t(s)>0$ such that $(t\star sh, t\star v)\in\cp_{a,b}^{-}$, where $t=t(s)$ is determined by
$$
\|\nabla v\|_2^2+\|\nabla h\|_2^2s^2=\frac{\gamma_p}{p}\mu_1\|h\|_p^ps^pt^{\gamma_p-2}+\frac{\gamma_q}{q}\mu_2\|v\|_q^qt^{\gamma_q-2}
+\gamma_{\alpha+\beta}\nu\int_{\R^N}|h|^{\alpha}|v|^{\beta}{\rm d}xs^{\alpha}t^{\gamma_{\alpha+\beta}-2}
$$
and
$$
\begin{aligned}
\|\nabla v\|_2^2+\|\nabla h\|_2^2s^2<&\frac{\mu_1}{p}\gamma_p(\gamma_p-1)\|h\|_p^ps^pt^{\gamma_p-2}+\frac{\mu_2}{q}\gamma_q(\gamma_q-1)\|v\|_q^qt^{\gamma_q-2}\\
&+\gamma_{\alpha+\beta}(\gamma_{\alpha+\beta}-1)\nu\int_{\R^N}|h|^{\alpha}|v|^{\beta}{\rm d}xs^{\beta}t^{\gamma_{\alpha+\beta}-2}.
\end{aligned}
$$
By the implicit function theorem, $t(s)\in C^{1}$ locally near $s=0$ and
$$
t'(s)=\frac{P_h(s)}{Q_h(s)},
$$
where
$$
P_h(s)=2\|\nabla h\|_2^2s-\gamma_p\mu_1\|h\|_p^ps^{p-1}t^{\gamma_p-2}-\gamma_{\alpha+\beta}\alpha\nu\int_{\R^N}|h|^{\alpha}|v|^{\beta}{\rm d}xs^{\alpha-1}t^{\gamma_{\alpha+\beta}-2},
$$
$$
\begin{aligned}
Q_h(s)=&\frac{\gamma_p(\gamma_p-2)}{p}\mu_1\|h\|_p^ps^pt^{\gamma_p-3}+\frac{\gamma_q(\gamma_q-2)}{q}\mu_2\|v\|_q^qt^{\gamma_q-3}\\
&+\gamma_{\alpha+\beta}(\gamma_{\alpha+\beta}-2)\nu\int_{\R^N}|h|^{\alpha}|v|^{\beta}{\rm d}xs^{\alpha}t^{\gamma_{\alpha+\beta}-3}.
\end{aligned}
$$
\textbf{Case $1<\alpha<2$:} Similar to \cite[Lemma 7.3]{Jeanjean2024a}, for $s>0$ small,
$$t(s)=1-M_h s^{\alpha}(1+o(1))$$
and for any $\tau>0$,
$$
t^{\tau}(s)=1-\tau M_hs^{\alpha}(1+o(1)),
$$
where
$$
M_h=\frac{\gamma_{\alpha+\beta}\nu\int_{\R^N}|h|^{\alpha}|v|^{\beta}{\rm d}x}{(\gamma_q-2)\|\nabla v\|_2^2}.
$$
Thus,
$$
\begin{aligned}
&I(t\star sh,t\star v)-I(0,v)\\
=&\frac{1}{2}\|\nabla v\|_2^2(t^2-1)-\frac{\mu_2}{q}\|v\|_q^q(t^{\gamma_q}-1)+\frac{1}{2}\|\nabla h\|_2^2s^2t^2-\frac{\mu_1}{p}\|h\|_p^ps^pt^{\gamma_p}
-\nu\int_{\R^N}|h|^{\alpha}|v|^{\beta}{\rm d}x s^{\alpha}t^{\gamma_{\alpha+\beta}}\\
=&-\|\nabla v\|_2^2M_h s^{\alpha}(1+o(1))+\frac{\mu_2}{q}\gamma_q\|v\|_q^qM_hs^{\alpha}(1+o(1))-\nu\int_{\R^N}|h|^{\alpha}|v|^{\beta}{\rm d}x s^{\alpha}(1+o(1))\\
=&-\nu\int_{\R^N}|h|^{\alpha}|v|^{\beta}{\rm d}x s^{\alpha}(1+o(1)).
\end{aligned}
$$
Hence, for sufficiently small $s>0$,
$$
Z_{a,b}\leq I(t\star sh,t\star v)<I(0,v)=m_q^{\mu_2}(b).
$$

\textbf{Case $\alpha=2$:} Repeating this process, we obtain
$$
I(t\star sh,t\star v)-I(0,v)=\left(\frac{1}{2}\|\nabla h\|_2^2-\nu\int_{\R^N}|h|^{2}|v|^{\beta}{\rm d}x\right)s^2(1+o(1)).
$$
By
$$\nu>\nu_{\mu_2,b,q,N,\beta}:=\frac{1}{2}\inf\limits_{h\in H^1(\R^N)\setminus\{0\}}\frac{\int_{\R^N}|\nabla h|^2{\rm d}x}{\int_{\R^N}|h|^2|v|^{\beta}{\rm d}x},$$
we can choose $h\in H^1(\R^N)\setminus\{0\}$ such that
$$
\frac{1}{2}\|\nabla h\|_2^2-\nu\int_{\R^N}|h|^{2}|v|^{\beta}{\rm d}x<0.
$$
Hence, for $\alpha=2$,
$$
Z_{a,b}\leq I(t\star sh,t\star v)<I(0,v)=m_q^{\mu_2}(b).
$$
\ep

\subsection{The minimizing problem $C_{a,b}=\inf\limits_{\cp_{a,b}^{+}}I(u,v)$}
\bl\label{le:240830-2}
Assume condition \eqref{condition:240901-1}, $N\geq1$ and $\widetilde T_{a,b}<c_0$. Let $m(a,b)$ be defined as in \eqref{def:240901-1}. If $C_{a,b}<m_p^{\mu_1}(a)$, then $m(a,b)$ is attained at some $(u,v)\in V_{R}$ satisfying $u \geq 0$ and $v \geq 0$ in $\mathbb{R}^N$,
with $u \not\equiv 0$ and $v \not\equiv 0$.
\el
\bp
Let $\{(u_n,v_n)\} \subset V_R$ be a minimizing sequence for $m(a,b)$. Consider the Schwarz rearrangement $(u^*_n,v^*_n)$, which by the properties of Schwarz rearrangement also forms a minimizing sequence. Without loss of generality, we may therefore assume that $u_n$ and $v_n$ are nonnegative, radially symmetric, non-increasing functions.

Applying Ekeland's variational principle, we obtain a Palais-Smale sequence satisfying $I(u_n,v_n) \to m(a,b)$ and $\|I'|_{V_R}(u_n,v_n)\| \to 0$. By Remark \ref{br:240729-1}, we have $(\|\nabla u_n\|_2^2+\|\nabla v_n\|_2^2)^{\frac{1}{2}}\leq R_0<R$ for large $n$. Consequently, there exists $\{(\lambda_{1,n}, \lambda_{2,n})\} \subset \mathbb{R}^2$ such that:
\beq\label{eq:240830-2}
\begin{cases}
-\Delta u_n+\lambda_{1,n}u_n=\mu_1|u_n|^{p-2}u_n+\nu\alpha |u_n|^{\alpha-2}u_n|v_n|^{\beta}+o_n(1) ~\hbox{in}~ H^{-1}(\R^N),\\
-\Delta v_n+\lambda_{2,n}v_n=\mu_2|v_n|^{q-2}v_n+\nu\beta |u_n|^{\alpha}|v_n|^{\beta-2}v_n+o_n(1) ~\hbox{in}~ H^{-1}(\R^N).
\end{cases}
\eeq
Since $\{(u_n,v_n)\}$ is bounded in $V_R$, there exists a subsequence converging weakly to $(u,v) \in V_R$ in $E_{rad}$. The radial symmetry implies $u_n \to u$ and $v_n \to v$ in $L^\eta(\mathbb{R}^N)$ for all $\eta \in (2, 2^*)$, $N \geq 1$. Applying the Br\'{e}zis-Lieb Lemma and H\"{o}lder's inequality, we have:
$$
\int_{\R^N}|u_n|^{\alpha}|v_n|^{\beta}{\rm d}x=\int_{\R^N}|u|^{\alpha}|v|^{\beta}{\rm d}x+\int_{\R^N}|u_n-u|^{\alpha}|v_n-v|^{\beta}{\rm d}x+o(1),
$$
and
$$
\int_{\R^N}|u_n-u|^{\alpha}|v_n-v|^{\beta}{\rm d}x\leq\|u_n-u\|^{\alpha+\beta}_{L^{\alpha+\beta}(\R^N)}\|v_n-v\|^{\alpha+\beta}_{L^{\alpha+\beta}(\R^N)}\to0,
$$
which implies that $\int_{\R^N}|u_n|^{\alpha}|v_n|^{\beta}{\rm d}x\to \int_{\R^N}|u|^{\alpha}|v|^{\beta}{\rm d}x$.

By the weak lower semicontinuity of the norm,
$$
m(a,b)=\liminf\limits_{n\to\infty}I(u_n,v_n)\geq I(u,v)\geq m(a,b).
$$
Thus, $I(u,v)=m(a,b)$. We now verify that $u\not\equiv0, v\not\equiv0$.

\noindent\textbf{Case 1: $u \equiv 0$, $v \equiv 0$:} This implies $m(a,b) = 0$, a contradiction.

\noindent\textbf{Case 2: $u \not\equiv 0$, $v \equiv 0$:} Since $m_p^{\mu_1}(a)=\inf\limits_{S_a}I(u,0)$, it follows that
$$m(a,b)=I(u,0)=\frac{1}{2}\|\nabla u\|_2^2-\frac{\mu_1}{p}\|u\|_p^p\geq m_p^{\mu_1}(\|u\|_2^2)\geq m_p^{\mu_1}(a),$$
contradicting $m(a,b)=C_{a,b}<m_p^{\mu_1}(a)$.

\noindent\textbf{Case 3: $u \equiv 0$, $v \not\equiv 0$:} By \eqref{eq:240830-2} and $v\not\equiv0$, we conclude that $\{\lambda_{2,n}\}$ is bounded. Up to a subsequence,  $\lambda_{2,n} \to \lambda_2$, and by a standard argument, $v$ satisfies
$$
-\Delta v+\lambda_2 v=\mu_2 v^{q-1} ~\hbox{in}~ \R^N.
$$
Thus $v\in\cp_{0,b}$. Then
$$m(a,b)=I(0,v)=\frac{1}{2}\|\nabla v\|_2^2-\frac{\mu_2}{q}\|v\|_q^q\geq m_q^{\mu_2}(\|v\|_2^2)\geq m_q^{\mu_2}(b)>0,$$
again a contradiction.
\ep

\br\label{br:240830-1}
In Lemma \ref{le:240830-2}, the conditions $u\not\equiv0$ and $v\not\equiv0$ imply that the sequences $\{\lambda_{i,n}\}$ ($i=1,2$) are bounded. Therefore, let $(u,v)$ be the minimizer of $I|_{V_R}$ corresponding to $ m(a,b)$, there exist $\lambda_1, \lambda_2\in \R$ satisfying
\beq\label{eq:240830-3}
\begin{cases}
-\Delta u+\lambda_1 u=\mu_1 u^{p-1}+\nu\alpha u^{\alpha-1}v^{\beta} ~\hbox{in}~ \R^N,\\
-\Delta v+\lambda_2 v=\mu_2 v^{q-1}+\nu\beta u^{\alpha}v^{\beta-1} ~\hbox{in}~ \R^N,\\
u\gneqq0, v\gneqq0.
\end{cases}
\eeq
By the Liouville-type result \cite[Lemma A.2]{Ikoma2014}, for dimensions $1\leq N\leq4$, we conclude that $\lambda_1>0$ and $\lambda_2>0$.
\er

\subsection{Proof of Theorem \ref{th:240829-1}}
To establish whether $(u,v) \in T(a,b)$, we adapt the argument of \cite[Lemma 8.2]{Jeanjean2024a} to obtain the following:
\bl\label{le:240830-3}
Let $(u,v)$ be given by Lemma \ref{le:240830-2} and $(\lambda_1,\lambda_2)$ by Remark \ref{br:240830-1}. Then:
\begin{itemize}
\item[(i)] If $\lambda_1 > 0$, then $u \in S_a$.
\item[(ii)] If $\lambda_2 > 0$, then $v \in S_b$.
\end{itemize}
\el
\bp
Suppose $\lambda_1>0$ and $\|u\|_2^2=\delta\in(0,a)$. For sufficiently small $s > 0$, we have $((1+s)u, v) \in D_a \times D_b$. By Lemma \ref{le:240824-2}, there exists a unique $t(s) > 0$ such that $(t(s) \star (1+s)u, t(s) \star v) \in \mathcal{P}_{a,b}^{+}$, where $t(s)$ satisfies:
\begin{align*}
(1+s)^2\|\nabla u\|_2^2+\|\nabla v\|_2^2=&\frac{\mu_1}{p}\gamma_p\|u\|_p^p(1+s)^{p}t^{\gamma_p-2}+\frac{\mu_2}{q}\gamma_q\|v\|_q^qt^{\gamma_q-2}\\
+&\gamma_{\alpha+\beta}\nu\int_{\R^N}u^{\alpha}v^{\beta}{\rm d}x(1+s)^{\alpha}t^{\gamma_{\alpha+\beta}-2}
\end{align*}
and
\begin{align*}
(1+s)^2\|\nabla u\|_2^2+\|\nabla v\|_2^2>&\frac{\mu_1}{p}\gamma_p(\gamma_p-1)\|u\|_p^p(1+s)^{p}t^{\gamma_p-2}+\frac{\mu_2}{q}\gamma_q(\gamma_q-1)\|v\|_q^qt^{\gamma_q-2}\\
+&\gamma_{\alpha+\beta}(\gamma_{\alpha+\beta}-1)\nu\int_{\R^N}u^{\alpha}v^{\beta}{\rm d}x(1+s)^{\alpha}t^{\gamma_{\alpha+\beta}-2}.
\end{align*}
By the implicit function theorem, $t(s)\in C^1$. Since
\begin{align*}
I(t\star(1+s)u,t\star v)=&\frac{1}{2}(\|\nabla u\|_2^2(1+s)^2+\|\nabla v\|_2^2)t^2-\frac{\mu_1}{p}\|u\|_p^p(1+s)^pt^{\gamma_p}-\frac{\mu_2}{q}\|v\|_q^qt^{\gamma_q}\\
-&\nu\int_{\R^N}u^{\alpha}v^{\beta}{\rm d}x(1+s)^{\alpha}t^{\gamma_{\alpha+\beta}},
\end{align*}
it follows that
\begin{align*}
\frac{d}{ds}I(t\star(1+s)u,t\star v)=&t^2(1+s)\|\nabla u\|_2^2+tt'(\|\nabla u\|_2^2(1+s)^2+\|\nabla v\|_2^2)\\
-&\mu_1\|u\|_p^pt^{\gamma_p}(1+s)^{p-1}-\frac{\mu_1}{p}\gamma_p\|u\|_p^p(1+s)^{p}t^{\gamma_p-1}t'-\frac{\mu_2}{q}\gamma_q\|v\|_q^qt^{\gamma_q-1}t'\\
-&\gamma_{\alpha+\beta}\nu\int_{\R^N}u^{\alpha}v^{\beta}{\rm d}x(1+s)^{\alpha}t^{\gamma_{\alpha+\beta}-1}t'-\alpha\nu\int_{\R^N}u^{\alpha}v^{\beta}{\rm d}x(1+s)^{\alpha-1}t^{\gamma_{\alpha+\beta}},
\end{align*}
where $t'=t'(s)$. Since $t(0)=1$ and $P(u,v)=0$, we obtain at $s=0$:
\begin{align*}
\frac{d}{ds}I(t\star(1+s)u,t\star v)|_{s=0}=&\|\nabla u\|_2^2-\mu_1\|u\|_p^p-\alpha\nu\int_{\R^N}u^{\alpha}v^{\beta}{\rm d}x\\
=&-\lambda_1\|u\|_2^2<0.
\end{align*}
Hence, for sufficiently small $s>0$,
$$
C_{a,b}\leq I(t\star(1+s)u,t\star v)<I(u,v)=m(a,b)=C_{a,b},
$$
a contradiction. Thus $u\in S_a$. An analogous argument shows $v\in S_b$ when $\lambda_2>0$.
\ep

\noindent\textbf{Proof of Theorem \ref{th:240829-1}:}
By the definition of $\widetilde{T}_{a,b}$, for any $b>0$ and $\nu>\nu_{\mu_1,p,a,N,\alpha}$, there exists $a(b, \nu)>0$ sufficiently small such that $\widetilde{T}_{a,b} < \tilde{c}_0 := \min\{\tilde{\alpha}_1,\tilde{\alpha}_2\}$ whenever $a < a(b, \nu)$.

For $\beta = 2$, Lemma \ref{le:240830-1}-(i) implies $C_{a,b} < m_p^{\mu_1}(a)$ when $\nu > \nu_{\mu_1,p,a,N,\alpha}$. By Lemma \ref{le:240830-2}, the minimization problem
$$\inf\limits_{V_R}I(u,v)=m(a,b)=C_{a,b}=\inf\limits_{\cp_{a,b}^{+}}I(u,v)$$
is attained at a pair positive, radially symmetric decreasing function $(u,v) \in V_R$. Combining Remark \ref{br:240830-1} and Lemma \ref{le:240830-3}, we conclude $(u,v) \in T(a,b)$. Thus, $(u,v)$ is a positive normalized ground state solution of the system \eqref{eq:240722-1}-\eqref{eq:240722-2}, and is also a local minimizer. Finally, for dimensions $N=1,2$, \cite[Lemma 7.1]{Jeanjean2024a} implies $\nu_{\mu_1,p,a,N,\alpha}=0$.
\hfill$\Box$

\subsection{Mountain pass geometric structure}
\bl\label{le:240830-4}
Under the hypotheses of Theorem \ref{th:240829-1}, there exists constants $k_0:=h(R)>0$ and $\rho<\min\{R,\sqrt{2k_0}\}$ such that
\beq\label{def:240830-2}
M(a,b):=\inf\limits_{\gamma\in\Gamma}\max\limits_{t\in[0,1]}I(\gamma(t))\geq k_0>\sup\limits_{\gamma\in\Gamma}\max\{I(\gamma(0)), I(\gamma(1))\},
\eeq
where $\Gamma$ denotes the set of continuous paths
$$
\begin{aligned}
\Gamma:=\{\gamma\in C([0,1],D_a\times D_b): &\gamma(t)=(\gamma_1(t),\gamma_2(t)),(\|\nabla \gamma_1(0)\|_2^2+\|\nabla \gamma_2(0)\|_2^2)^{\frac{1}{2}}<\rho,\\
&I(\gamma(1))<2m(a,b)\}.
\end{aligned}
$$
\el
\bp
First, we verify that $\Gamma\neq\emptyset$ and $M(a,b)$ is well-defined. Let $(u,v)$ be the local minimizer of $I|_{T(a,b)}$. By the fiber map \eqref{eq:240731-4}, $t \star (u,v) \in T(a,b)$ and
\begin{align*}
I(t\star(u,v))&=\frac{1}{2}(\|\nabla u\|_2^2+\|\nabla v\|_2^2)t^2-\frac{\mu_1}{p}\|u\|_p^pt^{\gamma_p}-\frac{\mu_2}{q}\|v\|_q^qt^{\gamma_q}-\nu\int_{\R^N}|u|^{\alpha}|v|^{\beta}{\rm d}xt^{\gamma_{\alpha+\beta}}\nonumber\\
&\to-\infty, ~\hbox{as}~ t\to+\infty.
\end{align*}
Thus, for sufficiently large $T>0$ large and small $t_0>0$, $I(T\star(u,v))<2m(a,b)<0$ and $t_0(\|\nabla u\|_2^2+\|\nabla v\|_2^2)^{\frac{1}{2}}<\rho$. Define the path $\gamma(t)=(t_0+(T-t_0)t)\star (u,v)$, which belongs to $\Gamma$.

For any $\gamma\in\Gamma$, we have
$$I(\gamma(0))\leq\frac{1}{2}(\|\nabla \gamma_1(0)\|_2^2+\|\nabla \gamma_2(0)\|_2^2)<\frac{1}{2}\rho^2<k_0.$$
Consequently,
$$
\sup\limits_{\gamma\in\Gamma}\max\{I(\gamma(0)), I(\gamma(1))\}<k_0.
$$
For any $\gamma\in\Gamma$, since $I(\gamma(1))<2m(a,b)$, it follows that $(\|\nabla \gamma_1(1)\|_2^2+\|\nabla \gamma_2(1)\|_2^2)^{\frac{1}{2}}>R$. Combining $(\|\nabla \gamma_1(0)\|_2^2+\|\nabla \gamma_2(0)\|_2^2)^{\frac{1}{2}}<\rho<R$, the Intermediate Value Theorem implies there exists $t_0\in(0,1)$ such that $(\|\nabla \gamma_1(t_0)\|_2^2+\|\nabla \gamma_2(t_0)\|_2^2)^{\frac{1}{2}}=R$.

By \eqref{eq:240901-1}, we have
$$
\max\limits_{t\in[0,1]}I(\gamma(t))\geq I(\gamma(t_0))\geq\inf\limits_{\partial V_R}I(u,v)\geq k_0>0.
$$
The arbitrariness of $\gamma$ implies $M(a,b)=\inf\limits_{\gamma\in\Gamma}\max\limits_{t\in[0,1]}I(\gamma(t))\geq k_0$.
\ep

\bl\label{le:240831-1}
Under the assumptions of Lemma \ref{le:240830-4}, there exists a nonnegative Palais-Smale-Pohozaev sequence $\{(u_n,v_n)\} \subset (D_a \times D_b) \cap E_{rad}$ for $I|_{D_a \times D_b}$ at level $M(a,b)$. That is,
\beq\label{eq:240830-4}
I(u_n,v_n)\to M(a,b), I|'_{D_a\times D_b}(u_n,v_n)\to0 ~\hbox{and}~ P(u_n,v_n)\to0.
\eeq
\el
\bp
The existence of such Palais-Smale-Pohozaev sequences under constrained conditions has been established in various settings \cite{Jeanjean1997,Bartsch2016}. Following these standard constructions in critical point theory, we omit the technical details.
\ep

\subsection{Estimation of the mountain pass level}
\bl\label{le:240830-7}
Under condition \eqref{condition:240901-1} and $\widetilde T_{a,b}<c_0$. let $Z_{a,b}$ and $M(a,b)$ be defined by \eqref{def:240830-1} and \eqref{def:240830-2} respectively. Then
$$
0<k_0\leq M(a,b)\leq Z_{a,b}.
$$
\el
\bp
For any $(u,v)\in \cp_{a,b}^{-}$, consider the fiber map $t\star (u,v)$. Select sufficiently small $t_0 > 0$ such that $t_0^2(\|\nabla u\|_2^2+\|\nabla v\|_2^2)<\rho^2$ and sufficiently large $t_1 > 0$ such that $I(t_1\star(u,v))<2m(a,b)$. Define the path
$$\gamma(t):=((1-t)t_0+tt_1)\star (u,v),\quad t\in[0,1].$$
By construction, $\gamma\in\Gamma$. Moreover, Lemma \ref{le:240824-2} implies
$$
M(a,b)\leq\max\limits_{t\in[0,1]}I(\gamma(t))=I(u,v).
$$
The arbitrariness of $(u,v)\in\cp_{a,b}^{-}$ yields $M(a,b)\leq Z_{a,b}$.
\ep

\subsection{Compactness and proof of Theorem \ref{th:240829-2}}
\bl\label{le:240830-5}
Let $\{(u_n,v_n)\}\subset E_{rad}$ be a nonnegative $(PSP)_{M(a,b)}$ sequence satisfying \eqref{eq:240830-4}. Then $\{(u_n,v_n)\}$ is bounded in $E$ and up to a subsequence, there exists $(u_0,v_0)\in E$ such that
$(u_n,v_n)\rightharpoonup (u_0,v_0)$ in $E_{rad}$.
\el
\bp
By \eqref{eq:240830-4}, we have
\beq\label{eq:240830-5}
I(u_n,v_n)=\frac{1}{2}(\|\nabla u_n\|_2^2+\|\nabla v_n\|_2^2)-\frac{\mu_1}{p}\|u_n\|_p^p-\frac{\mu_2}{q}\|v_n\|_q^q-\nu\int_{\R^N}u_n^{\alpha}v_n^{\beta}{\rm d}x=M(a,b)+o_n(1),
\eeq
\beq\label{eq:240830-6}
P(u_n,v_n)=\|\nabla u_n\|_2^2+\|\nabla v_n\|_2^2-\frac{\mu_1}{p}\gamma_p\|u_n\|_p^p-\frac{\mu_2}{q}\gamma_q\|v_n\|_q^q-\gamma_{\alpha+\beta}\nu\int_{\R^N}{u_n}^{\alpha}{v_n}^{\beta}{\rm d}x=o_n(1)
\eeq
and there exists $\{(\lambda_{1,n},\lambda_{2,n})\}$ such that system \eqref{eq:240830-2} holds.

For any $k\in\R$, combining \eqref{eq:240830-5} and \eqref{eq:240830-6}, we have
\begin{align}\label{eq:20250606-1}
&M(a,b)+o_n(1)=I(u_n,v_n)-kP(u_n,v_n)\nonumber\\
=&(\frac{1}{2}-k)(\|\nabla u_n\|_2^2+\|\nabla v_n\|_2^2)+\frac{\mu_1}{p}(k\gamma_p-1)\|u_n\|_p^p\nonumber\\
+&\frac{\mu_2}{q}(k\gamma_q-1)\|v_n\|_q^q+(k\gamma_{\alpha+\beta}-1)\nu\int_{\R^N}u_n^{\alpha}v_n^{\beta}{\rm d}x.
\end{align}

\noindent\textbf{Case $\alpha + \beta > 2 + \frac{4}{N}$.} Choose $k$ satisfying $\max\{\frac{1}{\gamma_q},\frac{1}{\gamma_{\alpha+\beta}}\}<k<\frac{1}{2}$, so that:
$$\frac{1}{2}-k, k\gamma_{\alpha+\beta}-1, k\gamma_q-1>0 ~\hbox{and}~ k\gamma_p-1<0.$$
From \eqref{eq:20250606-1}, we obtain
\beq\label{eq:240830-7}
\begin{aligned}
M(a,b)+o_n(1)>&(\frac{1}{2}-k)(\|\nabla u_n\|_2^2+\|\nabla v_n\|_2^2)+\frac{\mu_1}{p}(k\gamma_p-1)\|u_n\|_p^p\\
\geq&(\frac{1}{2}-k)(\|\nabla u_n\|_2^2+\|\nabla v_n\|_2^2)+\frac{\mu_1}{p}(k\gamma_p-1)\|u_n\|_2^{p-\gamma_p}(\|\nabla u_n\|_2^2+\|\nabla v_n\|_2^2)^{\frac{\gamma_p}{2}}\\
\geq&(\frac{1}{2}-k)(\|\nabla u_n\|_2^2+\|\nabla v_n\|_2^2)+\frac{\mu_1}{p}(k\gamma_p-1)a^{\frac{p-\gamma_p}{2}}(\|\nabla u_n\|_2^2+\|\nabla v_n\|_2^2)^{\frac{\gamma_p}{2}}.\\
\end{aligned}
\eeq
Since $\frac{\gamma_p}{2}<1$, it follows that $\|\nabla u_n\|_2^2+\|\nabla v_n\|_2^2<+\infty$.

\noindent\textbf{Case $\alpha + \beta < 2 + \frac{4}{N}$.} Select $\frac{1}{\gamma_q}<k<\frac{1}{2}$ so that:
$$\frac{1}{2}-k, k\gamma_q-1>0, k\gamma_{\alpha+\beta}-1 ~\hbox{and}~ k\gamma_p-1<0.$$
Set $r:=\alpha+\beta$. By H\"older and Gagliardo-Nirenberg inequalities:
\begin{align*}
\int_{\R^N}u_n^{\alpha}v_n^{\beta}{\rm d}x\leq&\|u_n\|_{r}^{\alpha}\|v_n\|_r^{\beta}\leq C_{N,r}\|u_n\|_2^{\frac{\alpha}{r}(r-\gamma_r)}\|\nabla u_n\|_2^{\frac{\alpha\gamma_r}{r}}\|v_n\|_2^{\frac{\beta}{r}(r-\gamma_r)}\|\nabla v_n\|_2^{\frac{\beta\gamma_r}{r}}\\
\leq&C_{N,r}a^{\frac{\alpha(r-\gamma_r)}{2r}}b^{\frac{\beta(r-\gamma_r)}{2r}}(\frac{\alpha}{r}\|\nabla u_n\|_2^2+\frac{\beta}{r}\|\nabla v_n\|_2^2)^{\frac{\gamma_r}{2}}.
\end{align*}
Then from \eqref{eq:20250606-1},
\beq\label{eq:240830-8}
\begin{aligned}
M(a,b)+o_n(1)>&(\frac{1}{2}-k)(\|\nabla u_n\|_2^2+\|\nabla v_n\|_2^2)+\frac{\mu_1}{p}(k\gamma_p-1)\|u_n\|_p^p\\
&+(k\gamma_{r}-1)\nu\int_{\R^N}|u_n|^{\alpha}|v_n|^{\beta}{\rm d}x\\
\geq&(\frac{1}{2}-k)(\|\nabla u_n\|_2^2+\|\nabla v_n\|_2^2)+\frac{\mu_1}{p}(k\gamma_p-1)a^{\frac{p-\gamma_p}{2}}(\|\nabla u_n\|_2^2+\|\nabla v_n\|_2^2)^{\frac{\gamma_p}{2}}\\
&+(k\gamma_{r}-1)\nu D_3(\|\nabla u_n\|_2^2+\|\nabla v_n\|_2^2)^{\frac{\gamma_{r}}{2}},
\end{aligned}
\eeq
where $D_3=(\frac{\max\{\alpha,\beta\}}{r})^{\frac{\gamma_r}{2}}C_{N,r}a^{\frac{\alpha(r-\gamma_r)}{2r}}b^{\frac{\beta(r-\gamma_r)}{2r}}$.
Since $\frac{\gamma_p}{2}<1$ and $\frac{\gamma_{r}}{2}<1$, it follows that $\|\nabla u_n\|_2^2+\|\nabla v_n\|_2^2<+\infty$.

\noindent\textbf{Case $\alpha + \beta = 2 + \frac{4}{N}$.} Take $\frac{1}{\gamma_q}<k<\frac{1}{2}$. Then $\gamma_{\alpha+\beta}=2$ and
$$\frac{1}{2}-k>0, k\gamma_q-1>0, k\gamma_{\alpha+\beta}-1<0 ~\hbox{and}~ k\gamma_p-1<0.$$
We also obtain \eqref{eq:240830-8} and so
\begin{align*}
M(a,b)+o_n(1)\leq&(\frac{1}{2}-k)(1-2\nu D_3)(\|\nabla u_n\|_2^2+\|\nabla v_n\|_2^2)\\
&+\frac{\mu_1}{p}(k\gamma_p-1)a^{\frac{p-\gamma_p}{2}}(\|\nabla u_n\|_2^2+\|\nabla v_n\|_2^2)^{\frac{\gamma_p}{2}}.
\end{align*}
In this case, by choosing a smaller $c_0$ if needed, we can ensure $\nu D_3<\frac{1}{2}$. Since $\frac{\gamma_p}{2}<1$, it follows that $\|\nabla u_n\|_2^2+\|\nabla v_n\|_2^2<+\infty$.

In all cases, $\{(u_n,v_n)\}$ is bounded in $E$. Thus, up to a subsequence, there exists $(u_0,v_0) \in E_{rad}$ such that $(u_n,v_n) \rightharpoonup (u_0,v_0)$ in $E_{rad}$.
\ep

\bl\label{le:240830-6}
Let $N \geq 2$ and $\{(u_n,v_n)\} \subset E_{rad}$ be a nonnegative Palais-Smale-Pohozaev sequence at level $M(a,b)$ from Lemma \ref{le:240830-5}. Denote by $(u_0,v_0)$ its weak limit. If $M(a,b) < m_q^{\mu_2}(b)$, then
$$
(u_n,v_n)\to (u_0,v_0) ~\hbox{in}~ D_0^{1,2}(\R^N)\times D_0^{1,2}(\R^N).
$$
Furthermore, $u_0 \geq 0$ with $u_0 \not\equiv 0$ and $v_0 \geq 0$ with $v_0 \not\equiv 0$.
\el
\bp
The weak convergence $(u_n,v_n) \rightharpoonup (u_0,v_0)$ in $E_{rad}$ holds by Lemma \ref{le:240830-5}. The Rellich compactness theorem implies that for all $\eta \in (2,2^*)$,
$$
u_n \to u_0 \quad \text{in} \quad L^{\eta}(\mathbb{R}^N), \quad v_n \to v_0 \quad \text{in} \quad L^{\eta}(\mathbb{R}^N),
$$
and
$$\int_{\R^N}|u_n|^{\alpha}|v_n|^{\beta}{\rm d}x\to \int_{\R^N}|u_0|^{\alpha}|v_0|^{\beta}{\rm d}x.$$
Next we {\bf claim} $u_0\not\equiv0, v_0\not\equiv0$.

If $u_0\equiv0$ and $v_0\equiv0$, then \eqref{eq:240830-6} implies
\begin{align*}
\|\nabla u_n\|_2^2+\|\nabla v_n\|_2^2=&\frac{\mu_1}{p}\gamma_p\|u_n\|_p^p+\frac{\mu_2}{q}\gamma_q\|v_n\|_q^q+\gamma_{\alpha+\beta}\nu\int_{\R^N}{u_n}^{\alpha}{v_n}^{\beta}{\rm d}x+o_n(1)\\
=&o_n(1).
\end{align*}
Thus $I(u_n,v_n)=o_n(1)$, contradicting $k_0>0$ in \eqref{def:240830-2}.

If $u_0\not\equiv0$ and $v_0\equiv0$. By \eqref{eq:240830-2}, $\{\lambda_{1,n}\}$ is bounded. Up to a subsequence, $\lambda_{1,n}\to\bar \lambda_1$ and $u_0$ satisfies
$$
-\Delta u_0+\bar\lambda_1 u_0=\mu_1 u_0^{p-1}.
$$
From \eqref{eq:240830-6} and $P(u_0,0)=0$, we have
\begin{align*}
\|\nabla u_n\|_2^2+\|\nabla v_n\|_2^2=&\frac{\mu_1}{p}\gamma_p\|u_n\|_p^p+\frac{\mu_2}{q}\gamma_q\|v_n\|_q^q+\gamma_{\alpha+\beta}\nu\int_{\R^N}{u_n}^{\alpha}{v_n}^{\beta}{\rm d}x+o_n(1)\\
=&\frac{\mu_1}{p}\gamma_p\|u_0\|_p^p+o_n(1)\\
=&\|\nabla u_0\|_2^2+o_n(1).
\end{align*}
By the weak lower semicontinuity of the norm,
$$(u_n,v_n)\to (u_0,0) ~\hbox{in}~D_0^{1,2}(\R^N)\times D_0^{1,2}(\R^N).$$
Since $u_0\in\cp_{a,0}$, we obtain
$$
M(a,b)=I(u_n,v_n)+o_n(1)=\frac{1}{2}\|\nabla u_0\|_2^2-\frac{\mu_1}{p}\|u_0\|_p^p+o_n(1)=m_p^{\mu_1}(\|u_0\|_2^2)+o_n(1),
$$
a contradiction.

If $u_0\equiv0$ and $v_0\not\equiv0$. \eqref{eq:240830-2} implies $\{\lambda_{2,n}\}$ is bounded. Up to a subsequence, $\lambda_{2,n}\to\bar \lambda_2$. And $v_0$ satisfies
$$
-\Delta v_0+\bar\lambda_2 v_0=\mu_2 u_0^{q-1} ~\text{in}~ \R^N.
$$
Analogously to the case $u_0\not\equiv0,v_0\equiv0$,
$$(u_n,v_n)\to (0,v_0) ~\hbox{in}~D_0^{1,2}(\R^N)\times D_0^{1,2}(\R^N).$$
Since $v_0\in\cp_{0,b}$, we have
$$
M(a,b)=I(u_n,v_n)+o_n(1)=\frac{1}{2}\|\nabla v_0\|_2^2-\frac{\mu_2}{q}\|v_0\|_p^p+o_n(1)=m_q^{\mu_2}(\|v_0\|_2^2)+o_n(1)\geq m_q^{\mu_2}(b)+o_n(1),
$$
contradicting $M(a,b)< m_q^{\mu_2}(b)$.

Thus $u_0\not\equiv0$ and $v_0\not\equiv0$. By \eqref{eq:240830-6}, $\{(\lambda_{1,n}, \lambda_{2,n})\}$ is bounded. Up to a subsequence, there exists $\{(\bar \lambda_{1},\bar \lambda_{2})\}$ such that
\beq\label{def:240830-3}
(\lambda_{1,n},\lambda_{2,n})\to (\bar \lambda_{1},\bar \lambda_{2}) ~\hbox{in}~ \R^2.
\eeq
Hence $(u_0,v_0,\bar\lambda_1,\bar\lambda_2)$ satisfies \eqref{eq:240830-3} and
\beq\label{eq:240830-9}
P(u_0,v_0)=\|\nabla u_0\|_2^2+\|\nabla v_0\|_2^2-\frac{\mu_1}{p}\gamma_p\|u_0\|_p^p-\frac{\mu_2}{q}\gamma_q\|v_0\|_q^q-\gamma_{\alpha+\beta}\nu\int_{\R^N}{u_0}^{\alpha}{v_0}^{\beta}{\rm d}x=0.
\eeq
Combining \eqref{eq:240830-6} and \eqref{eq:240830-9}, we conclude $(u_n,v_n)\to (u_0,v_0)$ in $D_0^{1,2}(\R^N)\times D_0^{1,2}(\R^N)$.
\ep
\bl\label{le:240831-2}
Let $(u_0,v_0)$ and $(\bar\lambda_1,\bar\lambda_2)$ be given by Lemma \ref{le:240830-5} and \eqref{def:240830-3}. Then:
\begin{itemize}
\item [(i)] If $\bar\lambda_1>0$, then $u\in S_a$.
\item [(ii)] If $\bar\lambda_2>0$, then $v\in S_b$.
\end{itemize}
\el
\bp
By Lemmas \ref{le:240830-7} and \ref{le:240830-6}, we have $0<k_0\leq I(u_0,v_0)=M(a,b)\leq Z_{a,b}$. Thus, Lemma \ref{le:240824-2} implies $(u_0,v_0)\in\cp_{a,b}^{-}$.

Assume $\bar\lambda_1>0$ and $\delta:=\|u_0\|_2^2\in(0,a)$. For sufficiently small $s>0$, $((1+s)u_0,v_0)\in D_a\times D_b$. By Lemma \ref{le:240824-2}, there exists a unique $t=t(s)>0$ such that $(t\star(1+s)u_0,t\star v_0)\in\cp_{a,b}^{-}$.
Following the argument in Lemma \ref{le:240830-3}, we obtain
$$
\begin{aligned}
\frac{d}{ds}I(t\star(1+s)u_0,t\star v_0)|_{s=0}=&\|\nabla u_0\|_2^2-\mu_1\|u_0\|_p^p-\alpha\nu\int_{\R^N}u_0^{\alpha}v_0^{\beta}{\rm d}x\\
=&-\bar\lambda_1\|u_0\|_2^2<0.
\end{aligned}
$$
Thus, for sufficiently small $s>0$, we have
$$
Z_{a,b}\leq I(t\star(1+s)u_0,t\star v_0)<I(u_0,v_0)=M(a,b)\leq Z_{a,b},
$$
a contradiction. Therefore $u\in S_a$. The case $\bar{\lambda}_2 > 0$ similarly implies $v_0\in S_b$.
\ep

\noindent\textbf{Proof of Theorem \ref{th:240829-2}:}
Under the assumptions of Theorem \ref{th:240829-1}, let $(u,v)$ denote a normalized ground state. By Lemma \ref{le:240831-1}, there exists a nonnegative Palais-Smale-Pohozaev sequence $\{(u_n,v_n)\}$ at level $M(a,b)$ for $I|_{D_a \times D_b}$.

When either $1 < \alpha < 2$ or $\alpha = 2$ with $\nu > \nu_{\mu_2,q,b,N,\beta}$, Lemmas \ref{le:240830-1} and \ref{le:240830-7} yield the mountain pass level estimate:
$$
0 < M(a,b) \leq Z_{a,b} < m_{q}^{\mu_2}(b).
$$
Applying Lemmas \ref{le:240830-5} and \ref{le:240830-6}, we obtain strong convergence:
\begin{align*}
&(u_n, v_n) \to (u_0, v_0) ~\text{in}~ D^{1,2}_0(\mathbb{R}^N) \times D^{1,2}_0(\mathbb{R}^N),\\
&(\lambda_{1,n}, \lambda_{2,n}) \to (\bar{\lambda}_1, \bar{\lambda}_2) ~\text{in}~ \mathbb{R}^2
\end{align*}
with $u_0\gneqq0$ and $v_0\gneqq0$.

By the Liouville-type theorem \cite[Lemma A.2]{Ikoma2014}, $\bar{\lambda}_1 > 0$ and $\bar{\lambda}_2 > 0$. Lemma \ref{le:240831-2} then implies $(u_0, v_0) \in T(a,b)$, confirming it as a positive normalized solution of the system \eqref{eq:240722-1}-\eqref{eq:240722-2}. When $N=2$, \cite[Lemma 7.1]{Jeanjean2024a} implies $\nu_{\mu_1,p,a,N,\alpha}=0$ and $\nu_{\mu_2,q,b,N,\beta}=0$.
\hfill$\Box$

\s{The Sobolev critical case}\label{sec:240829-1}
\subsection{Minimization problem $C_{a,b}=\inf\limits_{\cp_{a,b}^{+}}I(u,v)$.}
To establish the existence of critical points attaining $C_{a,b}$, we first estimate this energy level.

\bl\label{le:240825-1}
Assume condition \eqref{condition:240722-1}, $N\geq3$ and $T(a,b)<\min\{\alpha_1,\alpha_2\}$, where $\alpha_1$, $\alpha_2$ are given by Lemmas \ref{le:240724-1} and \ref{le:240829-1}. Let $m_{p}^{\mu_1}(a)$ be as defined in Lemma \ref{le:240729-1}. Then:
\begin{itemize}
\item [(i)] If $1<\beta<2$, then $C_{a,b}<m_{p}^{\mu_1}(a)$.
\item [(ii)] If $\beta=2$, then $C_{a,b}<m_{p}^{\mu_1}(a)$ when $\nu>\nu_{\mu_1,p,a,N,\alpha}$, with $\nu_{\mu_1,p,a,N,\alpha}$ defined in \eqref{eq:240829-2}.
\end{itemize}
\el
\bp
Similar results can be found in \cite{Jeanjean2024a}. Let $(u,\lambda_0)\in S_a\times\R^+$ be the unique solution to \eqref{eq:240726-1} with $\mu=\mu_1$ and $\eta=p$. For any $h \in H^1(\mathbb{R}^N)$ satisfying $\|h\|_2^2 = 1$, we have $(u, sh) \in D_a \times D_b$ when $|s| \leq \sqrt{b}$. By Lemma \ref{le:240824-2}-(i), for sufficiently small $s \geq 0$, there exists a unique $t(s) > 0$ such that $(t(s) \star u, t(s) \star sh) \in \mathcal{P}_{a,b}^{+}$, determined by:
\begin{align*}
\|\nabla u\|_2^2+\|\nabla h\|_2^2s^2=\frac{\gamma_p}{p}\mu_1\|u\|_p^pt^{\gamma_p-2}+\frac{\gamma_q}{q}\mu_2\|h\|_q^qs^qt^{\gamma_q-2}+2^*\nu\int_{\R^N}|u|^{\alpha}|h|^{\beta}{\rm d}xs^{\beta}t^{2^*-2}
\end{align*}
and
\begin{align*}
\|\nabla u\|_2^2+\|\nabla h\|_2^2s^2>&\frac{\mu_1}{p}\gamma_p(\gamma_p-1)\|u\|_p^pt^{\gamma_p-2}+\frac{\mu_2}{q}\gamma_q(\gamma_q-1)\|h\|_q^qs^qt^{\gamma_q-2}\\
&+2^*(2^*-1)\nu\int_{\R^N}|u|^{\alpha}|h|^{\beta}{\rm d}xs^{\beta}t^{2^*-2}.
\end{align*}
The implicit function theorem guarantees $t(s)\in C^{1}$ locally near $s=0$ with
$$
t'(s)=\frac{P_h(s)}{Q_h(s)},
$$
where
\begin{align*}
P_h(s)&=2\|\nabla h\|_2^2s-\gamma_q\mu_2\|h\|_q^qs^{q-1}t^{\gamma_q-2}-2^*\beta\nu\int_{\R^N}|u|^{\alpha}|h|^{\beta}{\rm d}xs^{\beta-1}t^{2^*-2},\\
Q_h(s)&=\frac{\gamma_p(\gamma_p-2)}{p}\mu_1\|u\|_p^pt^{\gamma_p-3}+\frac{\gamma_q(\gamma_q-2)}{q}\mu_2\|h\|_q^qs^qt^{\gamma_q-3}+2^*(2^*-2)\nu\int_{\R^N}|u|^{\alpha}|h|^{\beta}{\rm d}xs^{\beta}t^{2^*-3}.
\end{align*}
\noindent\textbf{Case $1 < \beta < 2$:} For small $s>0$, $t(s)=1+o(1)$ and
\begin{align*}
&P_h(s)=-2^*\nu\beta\int_{\R^N}|u|^{\alpha}|h|^{\beta}{\rm d}xs^{\beta-1}(1+o(1)),\\
&Q_h(s)=\frac{\gamma_p(\gamma_p-2)}{p}\mu_1\|u\|_p^p(1+o(1)).
\end{align*}
Thus,
$$t'(s)=\frac{P_h(s)}{Q_h(s)}=-M_h\beta s^{\beta-1}(1+o(1)),$$
where
$$
M_h=\frac{2^*p\nu\int_{\R^N}|u|^{\alpha}|h|^{\beta}{\rm d}x}{\gamma_p(\gamma_p-2)\mu_1\|u\|_p^p}=\frac{2^*\nu\int_{\R^N}|u|^{\alpha}|h|^{\beta}{\rm d}x}{(\gamma_p-2)\|\nabla u\|_2^2}.
$$
Consequently, for sufficiently small $s>0$,
$$t(s)=1-M_h s^{\beta}(1+o(1)),$$
and for any $\tau>0$,
$$
t^{\tau}(s)=1-\tau M_hs^{\beta}(1+o(1)).
$$
The energy difference satisfies
$$
\begin{aligned}
&I(t\star u,t\star sh)-I(u,0)\\
=&\frac{1}{2}\|\nabla u\|_2^2(t^2-1)-\frac{\mu_1}{p}\|u\|_p^p(t^{\gamma_p}-1)+\frac{1}{2}\|\nabla h\|_2^2s^2t^2-\frac{\mu_2}{q}\|h\|_q^qs^qt^{\gamma_q}
-\nu\int_{\R^N}|u|^{\alpha}|h|^{\beta}{\rm d}x s^{\beta}t^{2^*}\\
=&-\|\nabla u\|_2^2M_h s^{\beta}(1+o(1))+\frac{\mu_1}{p}\gamma_p\|u\|_p^pM_hs^{\beta}(1+o(1))-\nu\int_{\R^N}|u|^{\alpha}|h|^{\beta}{\rm d}x s^{\beta}(1+o(1))\\
=&-\nu\int_{\R^N}|u|^{\alpha}|h|^{\beta}{\rm d}x s^{\beta}(1+o(1)).
\end{aligned}
$$
Thus for sufficiently small $s>0$, it follows that
$$
C_{a,b}\leq I(t\star u,t\star sh)<I(u,0)=m_p^{\mu_1}(a).
$$

\noindent\textbf{Case $\beta = 2$:} In this case, for small $s>0$ , we have
\begin{align*}
&P_h(s)=(2\|\nabla h\|_2^2-2 2^*\nu\int_{\R^N}|u|^{\alpha}|h|^2{\rm d}x)s(1+o(1)),\\
&Q_h(s)=\frac{\gamma_p(\gamma_p-2)}{p}\mu_1\|u\|_p^p(1+o(1))=(\gamma_p-2)\|\nabla u\|_2^2(1+o(1)).
\end{align*}
Then
$$
t'(s)=\overline M_h s(1+o(1)),
$$
where
$$
\overline M_h=:\frac{2\|\nabla h\|_2^2-2 2^*\nu\int_{\R^N}|u|^{\alpha}|h|^2{\rm d}x}{(\gamma_p-2)\|\nabla u\|_2^2}.
$$
Consequently, for small $s>0$,
$$
t(s)=1+\frac{1}{2}\overline M_h s^2(1+o(1)),
$$
and for any $\tau>0$,
$$t^{\tau}(s)=1+\frac{\tau}{2}\overline M_h s^2(1+o(1)).$$
The energy difference satisfies
\begin{align*}
&I(t\star u,t\star sh)-I(u,0)\\
=&\frac{1}{2}\|\nabla u\|_2^2(t^2-1)-\frac{\mu_1}{p}\|u\|_p^p(t^{\gamma_p}-1)+\frac{1}{2}\|\nabla h\|_2^2s^2t^2-\frac{\mu_2}{q}\|h\|_q^qs^qt^{\gamma_q}
-\nu\int_{\R^N}|u|^{\alpha}|h|^{2}{\rm d}x s^{2}t^{2^*}\\
=&\frac{1}{2}\|\nabla u\|_2^2\overline M_hs^2(1+o(1))-\frac{\mu_1}{2p}\gamma_p\|u\|_p^p\overline M_hs^2(1+o(1))+\frac{1}{2}\|\nabla h\|_2^2s^2(1+o(1))\\
&-\nu\int_{\R^N}|u|^{\alpha}|h|^{2}{\rm d}x s^{2}(1+o(1))\\
=&\left(\frac{1}{2}\|\nabla h\|_2^2-\nu\int_{\R^N}|u|^{\alpha}|h|^{2}{\rm d}x\right)s^2(1+o(1)).
\end{align*}
Since $\nu>\nu_{\mu_1,a,p,N,\alpha}:=\frac{1}{2}\inf\limits_{h\in H^1(\R^N)\setminus\{0\}}\frac{\int_{\R^N}|\nabla h|^2{\rm d}x}{\int_{\R^N}|u|^{\alpha}|h|^2{\rm d}x}$, we may select $h$ such that
$$
\frac{1}{2}\|\nabla h\|_2^2-\nu\int_{\R^N}|u|^{\alpha}|h|^{2}{\rm d}x<0.
$$
Hence for $\beta=2$,
$$
C_{a,b}\leq I(t\star u,t\star sh)<I(u,0)=m_p^{\mu_1}(a).
$$
\ep

\bl\label{le:240729-3}
Let $N \geq 3$, assume condition \eqref{condition:240722-1} and $T_{a,b} < \min\{\alpha_1, \alpha_2\}$. If $C_{a,b} < m_p^{\mu_1}(a)$, then the infimum $m(a,b)$ is attained at some $(u,v) \in V_R$ with $u \geq 0$, $u \not\equiv 0$ and $v \geq 0$, $v \not\equiv 0$.
\el
\bp
Let $\{(u_n,v_n)\} \subset V_R$ be a minimizing sequence for $m(a,b)$. Applying Schwarz rearrangement to each pair $(u_n, v_n)$, we obtain $(u_n^*, v_n^*)$. The properties of Schwarz rearrangement ensure that $\{(u_n^*, v_n^*)\}$ is also a minimizing sequence. Without loss of generality, we may thus assume $u_n$ and $v_n$ are nonnegative, radially symmetric, non-increasing functions.

By Ekeland's variational principle, there exists a Palais-Smale sequence $\{(u_n,v_n)\}$ satisfying:
$$
I(u_n,v_n) \to m(a,b) \quad \text{and} \quad \|I'|_{V_R}(u_n,v_n)\| \to 0.
$$

Remark \ref{br:240729-1} implies that for sufficiently large $n$, $(\|\nabla u_n\|_2^2+\|\nabla v_n\|_2^2)^{\frac{1}{2}}\leq R_0<R$. Consequently, there exists
$\{(\lambda_{1,n},\lambda_{2,n})\}\subset \R^2$ such that
\beq\label{eq:240729-1}
\begin{cases}
-\Delta u_n+\lambda_{1,n}u_n=\mu_1|u_n|^{p-2}u_n+\nu\alpha |u_n|^{\alpha-2}u_n|v_n|^{\beta}+o_n(1) ~\hbox{in}~ H^{-1}(\R^N),\\
-\Delta v_n+\lambda_{2,n}v_n=\mu_2|v_n|^{q-2}v_n+\nu\beta |u_n|^{\alpha}|v_n|^{\beta-2}v_n+o_n(1) ~\hbox{in}~ H^{-1}(\R^N).
\end{cases}
\eeq

Since $\{(u_n,v_n)\}\subset  V_R$ is bounded in $E_{rad}$, up to a subsequence, there exists $(u,v)\in V_R$ such that $(u_n,v_n)\rightharpoonup (u,v)$ in $E_{rad}$. Rellich compactness theorem implies that
$$u_n\to u ~\text{in}~ L^{\eta}(\R^N) ~\text{and}~ v_n\to v ~\text{in}~ L^{\eta}(\R^N), ~\text{for all}~ N\geq2 ~\text{and}~ \eta\in(2,2^*).$$
Suppose
\beq\label{eq:240801-5}
\begin{cases}
|\nabla (u_n-u)|^2\rightharpoonup \xi_1 ~\hbox{in}~ \cm(\R^N),\\
|\nabla (v_n-v)|^2\rightharpoonup \xi_2 ~\hbox{in}~ \cm(\R^N),\\
|u_n-u|^{2^*}\rightharpoonup \eta_1 ~\hbox{in}~ \cm(\R^N),\\
|v_n-v|^{2^*}\rightharpoonup \eta_2 ~\hbox{in}~ \cm(\R^N),\\
|u_n-u|^{\alpha}|v_n-v|^{\beta}\rightharpoonup \eta_3 ~\hbox{in}~ \cm(\R^N),
\end{cases}
\eeq
and define the concentration at infinity:
\beq\label{eq:240801-6}
\begin{cases}
\xi_{\infty,1}:=\lim\limits_{R\to\infty}\limsup\limits_{n\to\infty}\int_{|x|\geq R}|\nabla u_n|^2{\rm d}x\geq0,\\
\xi_{\infty,2}:=\lim\limits_{R\to\infty}\limsup\limits_{n\to\infty}\int_{|x|\geq R}|\nabla v_n|^2{\rm d}x\geq0,\\
\eta_{\infty,1}:=\lim\limits_{R\to\infty}\limsup\limits_{n\to\infty}\int_{|x|\geq R}|u_n|^{2^*}{\rm d}x\geq0,\\
\eta_{\infty,2}:=\lim\limits_{R\to\infty}\limsup\limits_{n\to\infty}\int_{|x|\geq R}|v_n|^{2^*}{\rm d}x\geq0,\\
\eta_{\infty,3}:=\lim\limits_{R\to\infty}\limsup\limits_{n\to\infty}\int_{|x|\geq R}|u_n|^{\alpha}|v_n|^{\beta}{\rm d}x\geq0.
\end{cases}
\eeq
By radial symmetry and the Strauss lemma,
$$
\eta_{\infty,1}=\eta_{\infty,2}=\eta_{\infty,3}=0.
$$
Furthermore, the Sobolev and H\"older inequalities yield:
\beq\label{eq:240801-9}
\|\eta_i\|^{\frac{2}{2^*}}\leq S^{-1}\|\xi_i\|,\quad i=1,2,
\eeq
\beq\label{eq:240827-1}
\|\eta_3\|\leq\|\eta_1\|^{\frac{\alpha}{2^*}}\|\eta_2\|^{\frac{\beta}{2^*}}.
\eeq

Applying the Br\'{e}zis-Lieb Lemma and Young inequality:
\beq\label{eq:240825-2}
\begin{aligned}
m(a,b)=&I(u_n,v_n)+o(1)\\
=&I(u,v)+\frac{1}{2}\left(\|\nabla (u_n-u)\|_2^2+\|\nabla (v_n-v)\|_2^2\right)-\nu\int_{\R^N}|u_n-u|^{\alpha}|v_n-v|^{\beta}{\rm d}x\\
\geq& I(u,v)+\frac{1}{2}\left(\|\nabla (u_n-u)\|_2^2+\|\nabla (v_n-v)\|_2^2\right)-\nu\left(\frac{\alpha}{2^*}\|u_n-u\|_{2^*}^{2^*}+\frac{\beta}{2^*}\|v_n-v\|_{2^*}^{2^*}\right)\\
\geq& I(u,v)+\left(\frac{1}{2}\|\nabla (u_n-u)\|_2^2-\frac{\nu\alpha}{2^*}S^{-\frac{2^*}{2}}\|\nabla (u_n-u)\|_2^{2^*}\right)\\
+&\left(\frac{1}{2}\|\nabla (v_n-v)\|_2^2-\frac{\nu\beta}{2^*}S^{-\frac{2^*}{2}}\|\nabla (v_n-v)\|_2^{2^*}\right).
\end{aligned}
\eeq

Note that
\begin{align*}
\|\nabla (u_n-u)\|_2^2=\|\nabla u_n\|_2^2-\|\nabla u\|_2^2+o(1)\leq R_0^2,\\
\|\nabla (v_n-v)\|_2^2=\|\nabla v_n\|_2^2-\|\nabla v\|_2^2+o(1)\leq R_0^2.
\end{align*}
By \eqref{eq:240825-1}, we have
$$
\begin{aligned}
&\frac{1}{2}\|\nabla (u_n-u)\|_2^2-\frac{\nu\alpha}{2^*}S^{-\frac{2^*}{2}}\|\nabla (u_n-u)\|_2^{2^*}\\
=&\|\nabla (u_n-u)\|_2^2\left(\frac{1}{2}-\frac{\nu\alpha}{2^*}S^{-\frac{2^*}{2}}\|\nabla (u_n-u)\|_2^{2^*-2}\right)\\
\geq&\|\nabla (u_n-u)\|_2^2\left(\frac{1}{2}-\frac{\nu\alpha}{2^*}S^{-\frac{2^*}{2}}R_0^{2^*-2}\right)
\geq0.
\end{aligned}
$$
Similarly, we obtain
$$\frac{1}{2}\|\nabla (v_n-v)\|_2^2-\frac{\nu\beta}{2^*}S^{-\frac{2^*}{2}}\|\nabla (v_n-v)\|_2^{2^*}\geq0.$$
Hence \eqref{eq:240825-2} implies
$$
m(a,b)=I(u_n,v_n)+o(1)\geq I(u,v)\geq m(a,b),
$$
which establishes $I(u,v)=m(a,b)<0$ and $(u_n,v_n)\to (u,v)$ in $D_{0}^{1,2}(\R^N)\times D_{0}^{1,2}(\R^N)$.

To complete the proof, we verify $u \not\equiv 0$ and $v \not\equiv 0$:

\textbf{Case 1: $u \equiv 0$, $v \equiv 0$.} Then $I(u,v)=0$, contradicting $ m(a,b)<0$.

\textbf{Case 2: $u \equiv 0$, $v \not\equiv 0$.} Following the approach in Lemma \ref{le:240830-2}, we have
$$m(a,b)=I(0,v)\geq m_q^{\mu_2}(\|v\|_2^2)\geq m_q^{\mu_2}(b)>0,$$
contradicting $ m(a,b)<0$.

\textbf{Case 3: $u \not\equiv 0$, $v \equiv 0$.} Similarly,
$$ m(a,b)=I(u,0)\geq m_p^{\mu_1}(\|u\|_2^2)\geq m_p^{\mu_1}(a),$$
contradicting the hypothesis $m(a,b)=C_{a,b}<m_p^{\mu_1}(a)$.
\ep

\br\label{br:240825-1}
In Lemma \ref{le:240729-3}, the results $u \not\equiv 0$ and $v \not\equiv 0$ imply that the sequences $\{\lambda_{i,n}\}$ ($i=1,2$) are bounded. Consequently, for the minimizer $(u,v)$ of $I|_{V_R}$ at $m(a,b)$, there exists $(\lambda_1,\lambda_2)\in \mathbb{R}^2$ satisfying the elliptic system:
\beq\label{eq:240731-1}
\begin{cases}
-\Delta u+\lambda_1 u=\mu_1 u^{p-1}+\nu\alpha u^{\alpha-1}v^{\beta} ~\hbox{in}~ \R^N,\\
-\Delta v+\lambda_2 v=\mu_2 v^{q-1}+\nu\beta u^{\alpha}v^{\beta-1} ~\hbox{in}~ \R^N,\\
u\gneqq0, v\gneqq0.
\end{cases}
\eeq
By the Liouville-type theorem \cite[Lemma A.2]{Ikoma2014}, for dimensions $N = 3, 4$, we conclude $\lambda_1 > 0$ and $\lambda_2 > 0$.
\er

\subsection{Proof of Theorem \ref{th:240722-1}}
\begin{lemma}\label{le:240825-2}
Let $(u,v)$ and $(\lambda_1,\lambda_2)$ be given by Lemma \ref{le:240729-3} and Remark \ref{br:240825-1}. Then:
\begin{itemize}
\item [(i)] If $\lambda_1 > 0$, then $u \in S_a$
\item [(ii)] If $\lambda_2 > 0$, then $v \in S_b$
\end{itemize}
\end{lemma}
\bp
The argument parallels Lemma \ref{le:240830-3} with $\gamma_{\alpha+\beta}$ replaced by $2^*$.
\ep

\noindent\textbf{Proof of Theorem \ref{th:240722-1}:}
Set $c_0 = \min\{\alpha_1, \alpha_2\}$. For any $b > 0$ and $\nu>0$, there exists $a(b,\nu) > 0$ sufficiently small such that $T_{a,b} < c_0$ whenever $a < a(b,\nu)$. Hence, under the hypotheses of Theorem \ref{th:240722-1}, $C_{a,b} < m_p^{\mu_1}(a)$ holds.

By Lemma \ref{le:240729-3}, the infimum $m(a,b)$ is attained at some $(u,v) \in V_R$ that is positive, radially symmetric, decreasing, and satisfies \eqref{eq:240731-1} with some $\lambda_1, \lambda_2$. Remark \ref{br:240825-1} gives $\lambda_1 > 0, \lambda_2 > 0$, and Lemma \ref{le:240825-2} implies $(u,v) \in T(a,b)$. Hence, this minimizer lies on the constraint set
$$
\{(u,v)\in T(a,b):(\|\nabla u\|_2^2+\|\nabla v\|_2^2)^{\frac{1}{2}}\leq R_0\}.
$$
Moreover, Lemma \ref{le:240825-3} establishes $I(u,v)=m(a,b)=C_{a,b}=\inf\limits_{\cp_{a,b}}I(u,v)$. Thus, $(u,v)$ is indeed a normalized ground state solution.
\hfill$\Box$

\subsection{Mountain pass geometric structure}
\bl\label{le:240801-1}
Under the assumptions of Theorem \ref{th:240722-1}, let $(u,v)$ be a local minimizer of $I|_{T(a,b)}$. Then there exists $k_0>0$ such that
$$
M(a,b):=\inf\limits_{\gamma\in\Gamma}\max\limits_{t\in[0,1]}I(\gamma(t))\geq k_0>0>\max\{I(\gamma(0)), I(\gamma(1))\},
$$
where $\Gamma$ denotes the set of continuous paths
$$\Gamma:=\{\gamma\in C([0,1],T(a,b)):\gamma(t)=(\gamma_1(t),\gamma_2(t)),\gamma(0)=(u,v),I(\gamma(1))<2m(a,b)\}.$$
\el
\bp
We first verify $\Gamma\neq\emptyset$ and $M(a,b)$ is well-defined. By the fiber map \eqref{eq:240731-4}, $t \star (u,v) \in T(a,b)$ for all $t > 0$, and
\begin{align*}
I(t\star(u,v))&=\frac{1}{2}(\|\nabla u\|_2^2+\|\nabla v\|_2^2)t^2-\frac{\mu_1}{p}\|u\|_p^pt^{\gamma_p}-\frac{\mu_2}{q}\|v\|_q^qt^{\gamma_q}-\nu\int_{\R^N}|u|^{\alpha}|v|^{\beta}{\rm d}xt^{2^*}\nonumber\\
&\to-\infty, ~\hbox{as}~ t\to+\infty.
\end{align*}
Thus, there exists $T > 0$ sufficiently large such that $I(T\star(u,v))<2m(a,b)<0$. Define $\gamma(t)=(1+(T-1)t)\star (u,v)$, which belongs to $\Gamma$.

For any $\gamma \in \Gamma$, note that $\gamma(0) = (u,v)$ implies
$$
(\|\nabla \gamma_1(0)\|_2^2+\|\nabla \gamma_2(0)\|_2^2)^{\frac{1}{2}}<R,
$$
while $I(\gamma(1)) < 2m(a,b)$ gives
$$
(\|\nabla \gamma_1(1)\|_2^2+\|\nabla \gamma_2(1)\|_2^2)^{\frac{1}{2}}>R.
$$
By the Intermediate Value Theorem, there exists $t_0 \in (0,1)$ such that
$$
(\|\nabla \gamma_1(t_0)\|_2^2 + \|\nabla \gamma_2(t_0)\|_2^2)^{\frac{1}{2}} = R.
$$

Set $k_0:=h(R)>0$. By \eqref{eq:240724-3}, we have
$$
\max\limits_{t\in[0,1]}I(\gamma(t))\geq I(\gamma(t_0))\geq\inf\limits_{\partial V_R}I(u,v)\geq k_0>0.
$$
The arbitrariness of $\gamma$ implies that $M(a,b)=\inf\limits_{\gamma\in\Gamma}\max\limits_{t\in[0,1]}I(\gamma(t))\geq k_0$.
\ep

\subsection{Existence of nonnegative $(PSP)_{M(a,b)}$ sequence}
\bl\label{le:20240927-1}
Under the hypotheses of Lemma \ref{le:240801-1}, there exists a nonnegative Palais-Smale-Pohozaev sequence $\{(u_n,v_n)\} \subset T(a,b) \cap E_{rad}$ for the constrained functional $I|_{T(a,b)}$, satisfying:
\beq\label{eq:240801-1}
I(u_n,v_n)\to M(a,b), I|'_{T(a,b)}(u_n,v_n)\to0 ~\hbox{and}~ P(u_n,v_n)\to0.
\eeq
\el
\bp
This kind of results is standard, we omit the detail.
\ep

\subsection{Estimation of the mountain pass level}
For the Sobolev-critical system \eqref{eq:240722-1}, obtaining precise estimates of $M(a,b)$ is crucial for establishing the compactness of Palais-Smale-Pohozaev sequences $\{(u_n,v_n)\}$ in $E_{rad}$. This is particularly challenging for mass-mixed systems, where mountain pass level estimation often constitutes the primary obstacle to proving multiplicity results.

Recent some studies have addressed normalized solutions for Schr\"{o}dinger systems with Sobolev-critical exponents:
\begin{itemize}
\item \cite{Vicentiu2023} treats $p=q=2^*$ with $2 < \alpha+\beta < 2 + \frac{4}{N}$
\item \cite{He2024} considers $2 < p,q < 2 + \frac{4}{N}$ and $\alpha+\beta = 2^*$
\end{itemize}
Our work focuses on the parameter
$$
2 < p < 2 + \tfrac{4}{N} < q < 2^*, \quad \alpha + \beta = 2^*.
$$
The methodologies developed in \cite{He2024,Vicentiu2023} provide essential foundations for our subsequent energy estimates.

Let $A_N:=[N(N-2)]^{\frac{N-2}{4}}$ and define $U_n(x):=\Theta_n(|x|)\in H_{rad}^{1}(\R^N)$, where
$$
\Theta_n(r)=A_N\left\{\begin{array}{lcl}
\left(\frac{n}{1+n^{2}r^{2}}\right)^{\frac{N-2}{2}}, &0\leq r<1;\\
\left(\frac{n}{1+n^2}\right)^{\frac{N-2}{2}}(2-r), &1\leq r<2;\\
0, &r\geq 2.
\end{array}
\right.
$$
Through direct computation, we derive the following asymptotic estimates:
\beq\label{eq:240802-1}
\begin{aligned}
\|U_n\|_2^2=O(\frac{\xi(n)}{n^2}), \quad  n\to\infty,
\end{aligned}
\eeq
\beq\label{0427_2}
\xi(n):=\int_0^n\frac{s^{N-1}}{(1+s^2)^{N-2}}ds=
\begin{cases}
O(n),\quad\quad &\hbox{if}~N=3;\\
O\big(\ln(1+n^2)\big), \quad &\hbox{if}~N=4.
\end{cases}
\eeq
\beq\label{0427_3}
\begin{aligned}
\|\nabla U_n\|_2^2&={\int_{\R^N}|\nabla U_n|^2{\rm d}x}=S^{\frac{N}{2}}+O\left(\frac{1}{n^{N-2}}\right), \quad n\to \infty.
\end{aligned}
\eeq
\beq\label{0427_4}
\begin{aligned}
\|U_n\|^{2^*}_{2^*}=S^{\frac{N}{2}}+O\left(\frac{1}{n^N}\right), \quad n\to \infty.
\end{aligned}
\eeq
Let $(u,v)$ be the local minimizer from Theorem \ref{th:240722-1}. By elliptic regularity theory, $u$ and $v$ are positive, radially decreasing, $C^2$  functions. Thus, we also have the following decay estimates (cf. \cite{Vicentiu2023}):
\beq\label{eq:240802-2}
\int_{\R^N}u U_n{\rm d}x=O\left(\frac{1}{n^{\frac{N-2}{2}}}\right),\int_{\R^N}v U_n{\rm d}x=O\left(\frac{1}{n^{\frac{N-2}{2}}}\right) \quad n\to \infty.
\eeq
\beq
\int_{\R^N}u U_n^{2^*-1}{\rm d}x=O\left(\frac{1}{n^{\frac{N-2}{2}}}\right), \int_{\R^N}v U_n^{2^*-1}{\rm d}x=O\left(\frac{1}{n^{\frac{N-2}{2}}}\right)\quad n\to \infty.
\eeq
For $t\geq0$ and $n\in N$, we define
$$
\begin{cases}
\Phi_{n,t}:=\frac{\sqrt{a}}{\|u+tU_n\|_2}(u+tU_n),\\
\Psi_{n,t}:=\frac{\sqrt{b}}{\|v+\sqrt{\frac{\beta}{\alpha}}tU_n\|_2}(v+\sqrt{\frac{\beta}{\alpha}}tU_n).
\end{cases}
$$
For each $n\in N$, $\{(\Phi_{n,t},\Psi_{n,t}):t\geq0\}$ is a curve in $T(a,b)\cap E_{rad}$ with $(\Phi_{n,0},\Psi_{n,0})=(u,v)$.

\bl\cite[Lemma 3.6]{He2024}\label{le:240802-1}
For any $T>0$, the following results hold uniformly for $t\in[0,T]$.
\begin{itemize}
\item [(i)] $\|u+tU_n\|_2^2=a+o_n(1), \|v+\sqrt{\frac{\beta}{\alpha}}tU_n\|_2^2=b+o_n(1)$;
\item [(ii)] $\|u+tU_n\|_p^p=\|u\|_p^p+o_n(1), \|v+\sqrt{\frac{\beta}{\alpha}}tU_n\|_q^q=\|v\|_q^q+o_n(1)$;
\item [(iii)] $\int_{\R^N}|u+tU_n|^{\alpha}|v+\sqrt{\frac{\beta}{\alpha}}tU_n|^{\beta}{\rm d}x=
    \int_{\R^N}u^{\alpha}v^{\beta}{\rm d}x+(\frac{\beta}{\alpha})^{\frac{\beta}{2}}S^{\frac{N}{2}}t^{2^*}+o_n(1)$;
\item [(iv)] $\|\nabla (u+tU_n)\|_2^2=\|\nabla u\|_2^2+S^{\frac{N}{2}}t^2+o_n(1), \|\nabla (v+\sqrt{\frac{\beta}{\alpha}}tU_n)\|_2^2=\|\nabla v\|_2^2+\frac{\beta}{\alpha}S^{\frac{N}{2}}t^2+o_n(1)$;
\item [(v)] For any $\eta\in\R$, $\left(\frac{\sqrt{a}}{\|u+tU_n\|_2}\right)^{\eta}=1+o_n(1)$, $\left(\frac{\sqrt{b}}{\|v+\sqrt{\frac{\beta}{\alpha}}tU_n\|_2}\right)^{\eta}=1+o_n(1)$, $\frac{d}{dt}\left(\left(\frac{\sqrt{a}}{\|u+tU_n\|_2}\right)^{\eta}\right)=o_n(1)$,
    $\frac{d}{dt}\left(\left(\frac{\sqrt{b}}{\|v+\sqrt{\frac{\beta}{\alpha}}tU_n\|_2\|_2}\right)^{\eta}\right)=o_n(1)$.
\end{itemize}
\el
Let $H_n:\R^+\to\R$ be a map defined by
\beq\label{eq:240802-5}
H_n(t):=I(\Phi_{n,t},\Psi_{n,t}),\quad t\geq0.
\eeq

\bl\label{le:240802-2}
Let $H_n(t)$ be defined as above. Then there exist $T>0$ and $n_{T}>0$ such that for all $n\geq n_T$, $H_n(T)<2m(a,b)$. Define $\gamma_n (t):=(\Phi_{n,Tt}, \Psi_{n,Tt})$, then $\gamma_n\in\Gamma$. Moreover, there exists $0<t_n<T$ such that $H_n(t_n)=\max\limits_{t>0}H_n(t)$ and $0<\inf\limits_{n\geq n_T}t_n<\sup\limits_{n\geq n_T}t_n<+\infty$.
\el
\bp
Observe that for any $n\in N$, $H_n(0)=I(\Phi_{n,0},\Psi_{n,0})=I(u,v)$. Fix $T > 0$ sufficiently large (to be determined). By Lemma \ref{le:240802-1}, uniformly for $t \in [0,T]$,
$$
\begin{aligned}
H_n(t)=&\frac{1}{2}(1+o_n(1))(\|\nabla u\|_2^2+S^{\frac{N}{2}}t^2+o_n(1))+\frac{1}{2}(1+o_n(1))(\|\nabla v\|_2^2+\frac{\beta}{\alpha}S^{\frac{N}{2}}t^2+o_n(1))\\
&-\frac{\mu_1}{p}(1+o_n(1))(\|u\|_p^p+o_n(1))-\frac{\mu_2}{q}(1+o_n(1))(\|v\|_q^q+o_n(1))\\
&-\nu(1+o_n(1))\left(\int_{\R^N}u^{\alpha}v^{\beta}{\rm d}x+(\frac{\beta}{\alpha})^{\frac{\beta}{2}}S^{\frac{N}{2}}t^{2^*}+o_n(1)\right)\\
=&m(a,b)+S^{\frac{N}{2}}\left[\frac{2^*}{2\alpha}t^2-\nu(\frac{\beta}{\alpha})^{\frac{\beta}{2}}t^{2^*}\right]+o_n(1).
\end{aligned}
$$
Since $2^* > 2$, we can select $T > 0$ large enough such that for all sufficiently large $n \geq n_T$, $H_n(T)<2m(a,b)$. Thus $\gamma_n\in\Gamma$.
By Lemma \ref{le:240801-1}, there exists $t_n\in(0,T)$ such that
$$
H_n(t_n)=\max\limits_{t>0}H_n(t)=\max\limits_{t\in[0,1]}I(\gamma_n(t))\geq k_0>0.
$$

We now establish $\inf\limits_{n\geq n_{T}}t_n>0$ by contradiction. Suppose there exists a subsequence $\{n_k\}$ with $n_k \to \infty$ and $t_{n_k} \to 0$. Lemma \ref{le:240802-1} gives
$$
(\Phi_{n_k,t_{n_k}}, \Psi_{n_k,t_{n_k}}) \to (u,v) \quad \text{in} \quad E,
$$
and thus
$$H_{n_k}(t_{n_k})=I(\Phi_{n_k,t_{n_k}},\quad \Psi_{n_k,t_{n_k}})\to I(u,v)=m(a,b)<0,$$
contradicting $H_n(t_n)\geq k_0>0$ for $n\geq n_{T}$. Therefore,
$$0<\inf\limits_{n\geq n_T}t_n<\sup\limits_{n\geq n_T}t_n\leq T<+\infty.$$
\ep

\br
For $n\geq n_{T}$, Lemma \ref{le:240802-2} implies that up to a subsequence, $t_n\to t^*\in (0,T)$.  Since $H_n'(t_n) = 0$, the limiting equation
\begin{equation}
\frac{1}{\alpha}t^* - \nu \left(\frac{\beta}{\alpha}\right)^{\frac{\beta}{2}} (t^*)^{2^*-1} = 0
\end{equation}
is obtained as $n \to \infty$. This equation admits the solution
$$t^*=\nu^{-\frac{N-2}{4}}\alpha^{\frac{4-(N-2)\alpha}{8}}\beta^{-\frac{(N-2)\beta}{8}}.$$
For full computational details, we refer to \cite[Lemma 3.9]{He2024}.
\er

To facilitate subsequent arguments, we adopt the following notations and identities from \cite{He2024}. For $N \in \{3,4\}$, equations \eqref{eq:240802-1}, \eqref{0427_2} and \eqref{eq:240802-2} imply the asymptotic expansions:
\beq\label{eq:240802-3}
\frac{\sqrt{a}}{\|u+t_n U_n\|_2}=1-l_1\frac{1}{n^{\frac{N-2}{2}}}+o\left(\frac{1}{n^{\frac{N-2}{2}}}\right),
\eeq
\beq\label{eq:240802-4}
\frac{\sqrt{b}}{\|v+\sqrt{\frac{\beta}{\alpha}}t U_n\|_2}=1-l_2\frac{1}{n^{\frac{N-2}{2}}}+o\left(\frac{1}{n^{\frac{N-2}{2}}}\right).
\eeq
Consequently,
\beq\label{eq:240803-1}
\left(\frac{\sqrt{a}}{\|u+t_n U_n\|_2}\right)^p=1-pl_1\frac{1}{n^{\frac{N-2}{2}}}+o\left(\frac{1}{n^{\frac{N-2}{2}}}\right),
\eeq
\beq\label{eq:240803-2}
\left(\frac{\sqrt{b}}{\|v+\sqrt{\frac{\beta}{\alpha}}t U_n\|_2}\right)^q=1-q l_2\frac{1}{n^{\frac{N-2}{2}}}+o\left(\frac{1}{n^{\frac{N-2}{2}}}\right).
\eeq
\bl\cite[Lemma 3.10]{He2024}\label{le:240803-1}
Let $N\in\{3,4\}$ and $l_1,l_2$ be the numbers given by \eqref{eq:240802-3} and \eqref{eq:240802-4} respectively, then the following hold:
$$
\begin{cases}
\int_{\R^N}u U_n {\rm d}x t_n=al_1\frac{1}{n^{\frac{N-2}{2}}}+o(\frac{1}{n^{\frac{N-2}{2}}}),\\
\int_{\R^N}v U_n{\rm d}x\sqrt{\frac{\beta}{\alpha}}t_n=b l_2\frac{1}{n^{\frac{N-2}{2}}}+o(\frac{1}{n^{\frac{N-2}{2}}}).
\end{cases}
$$
\el

\bl\cite[Lemma 3.12]{He2024}
Let $\alpha>1,\beta>1,\alpha+\beta>3$ and $0<L_1\leq L_2<\infty$, then there exists some $A_1>0$ small and $A_2>0$ large such that
$$
\begin{aligned}
&(t_1+s)^{\alpha}(t_2+s)^{\beta}-t_1^{\alpha}t_2^{\beta}-s^{\alpha+\beta}-\alpha t_1^{\alpha-1}t_2^{\beta}s-\beta t_2^{\beta-1}t_1^{\alpha}s\\
\geq&A_1s^{\alpha+\beta-1}-A_2 s^2, \forall (t_1,t_2,s)\in [L_1,L_2]^2\times\R^+.
\end{aligned}
$$
\el
\begin{corollary}\cite[Corollary 3.13]{He2024}\label{corollary:240803-1}
Let $3\leq N\leq 5$, $\alpha>1,\beta>1,\alpha+\beta=2^*$ and $U_n$ be given as above. Then
$$
\begin{aligned}
&(u+tU_n)^{\alpha}(v+\sqrt{\frac{\beta}{\alpha}}tU_n)^{\beta}-u^{\alpha}v^{\beta}-(\sqrt{\frac{\beta}{\alpha}})^{\beta}t^{2^*}U_n^{2^*}-\alpha u^{\alpha-1}v^{\beta}t U_n-\beta v^{\beta-1}u^{\alpha}\sqrt{\frac{\beta}{\alpha}}tU_n\\
\geq&A_1 (\sqrt{\frac{\beta}{\alpha}})^{\beta}t^{2^*-1}U_n^{2^*-1}-A_2(\sqrt{\frac{\beta}{\alpha}})^{\beta}t^2 U_n^2
\end{aligned}
$$
for all $t\geq0$ and $x\in\R^N$.
\end{corollary}

\bl\label{le:240803-3}
Under the hypotheses of Theorem \ref{th:240722-1}, let $H_n(t)$ be defined in \eqref{eq:240802-5} and $t_n$ be given in Lemma \ref{le:240802-2}. Then for sufficiently large $n$,
$$
H_n(t_n)<m(a,b)+\frac{2}{N-2}\nu^{-\frac{N-2}{2}}\alpha^{-\frac{(N-2)\alpha}{4}}\beta^{-\frac{(N-2)\beta}{4}}S^{\frac{N}{2}}.
$$
\el
\bp
Since $(u,v)$ is a local minimizer of $I|_{T(a,b)}$ satisfying \eqref{eq:240722-1} with $\|u\|_2^2 = a$, $\|v\|_2^2 = b$, we have:
\begin{align*}
\int_{\R^N}\nabla u \nabla U_n {\rm d}x=-\lambda_1\int_{\R^N}u U_n {\rm d}x+\mu_1\int_{\R^N} u^{p-1}U_n {\rm d}x+\nu\alpha\int_{\R^N} u^{\alpha-1}v^{\beta}U_n {\rm d}x,\\
\int_{\R^N}\nabla v \nabla U_n {\rm d}x=-\lambda_2\int_{\R^N}v U_n {\rm d}x+\mu_1\int_{\R^N} v^{q-1}U_n {\rm d}x+\nu\beta\int_{\R^N} u^{\alpha}v^{\beta-1}U_n {\rm d}x,
\end{align*}
where $\lambda_1,\lambda_2$ are given in Theorem \ref{th:240722-1}. Applying equations \eqref{eq:240802-3}-\eqref{eq:240803-2}, Lemma \ref{le:240803-1}, and Corollary \ref{corollary:240803-1}, the method established in \cite[Proposition 3.17]{He2024} extends directly to dimensions $N \in \{3,4\}$, and we consequently omit the computational details.
\ep

\begin{corollary}\label{le:240801-2}
Let $N\in\{3,4\}$. Assume condition \eqref{condition:240722-1} and let $c_0$ be as defined in Theorem \ref{th:240722-1}. For any $a,b$ satisfying \eqref{eq:240723-3}, the mountain pass level satisfies
\beq\label{eq:240801-10}
M(a,b)<m(a,b)+\frac{2}{N-2}\nu^{-\frac{N-2}{2}}\alpha^{-\frac{(N-2)\alpha}{4}}\beta^{-\frac{(N-2)\beta}{4}}S^{\frac{N}{2}}.
\eeq
\end{corollary}
\bp
By Lemma \ref{le:240802-2}, we have
$$M(a,b)\leq\max\limits_{t\in[0,1]}I(\Phi_{n,Tt}, \Psi_{n,Tt})=H_n(t_n).$$
The conclusion then follows directly from Lemma \ref{le:240803-3}.
\ep

\subsection{Proof of Theorem \ref{th:240722-2}}
\bl\label{le:240801-1}
Let $\{(u_n,v_n)\}$ be a nonnegative Palais-Smale-Pohozaev sequence at level $M(a,b)$ satisfying \eqref{eq:240801-1}. Then $\{(u_n,v_n)\}$ is bounded in $E$, and up to a subsequence, there exist $(u_0,v_0)\in E$ and $(\bar\lambda_1,\bar\lambda_2)\in \mathbb{R}^2$ such that:
\begin{itemize}
\item [(i)] $(u_n,v_n) \rightharpoonup (u_0,v_0)$ in $E_{rad}$, with $u_n \to u_0$ and $v_n \to v_0$ in $L^{\eta}(\mathbb{R}^N)$ for all $\eta \in (2,2^*)$.
\item [(ii)] The associated Lagrange multipliers satisfy $(\lambda_{1,n},\lambda_{2,n}) \to (\bar\lambda_1,\bar\lambda_2)$.
\item [(iii)] $(u_0,v_0,\bar\lambda_1,\bar\lambda_2)$ solves system \eqref{eq:240722-1}.
\end{itemize}
\el
\bp
By \eqref{eq:240801-1}, we have
\beq\label{eq:240801-2}
I(u_n,v_n)=\frac{1}{2}(\|\nabla u_n\|_2^2+\|\nabla v_n\|_2^2)-\frac{\mu_1}{p}\|u_n\|_p^p-\frac{\mu_2}{q}\|v_n\|_q^q-\nu\int_{\R^N}u_n^{\alpha}v_n^{\beta}{\rm d}x=M(a,b)+o_n(1),
\eeq
\beq\label{eq:240801-3}
P(u_n,v_n)=\|\nabla u_n\|_2^2+\|\nabla v_n\|_2^2-\frac{\mu_1}{p}\gamma_p\|u_n\|_p^p-\frac{\mu_2}{q}\gamma_q\|v_n\|_q^q-2^*\nu\int_{\R^N}{u_n}^{\alpha}{v_n}^{\beta}{\rm d}x=o_n(1)
\eeq
and there exists $\{(\lambda_{1,n},\lambda_{2,n})\}$ such that system \eqref{eq:240729-1} holds.

\eqref{eq:240801-2} and \eqref{eq:240801-3} yield for any $k\in\R$,
\begin{align*}
M(a,b)+o_n(1)=&I(u_n,v_n)-kP(u_n,v_n)\\
=&(\frac{1}{2}-k)(\|\nabla u_n\|_2^2+\|\nabla v_n\|_2^2)+\frac{\mu_1}{p}(k\gamma_p-1)\|u_n\|_p^p\\
+&\frac{\mu_2}{q}(k\gamma_q-1)\|v_n\|_q^q+(2^*k-1)\nu\int_{\R^N}|u_n|^{\alpha}|v_n|^{\beta}{\rm d}x.
\end{align*}
Select $k$ satisfying $\frac{1}{\gamma_q} < k < \frac{1}{2}$, ensuring:
$$\frac{1}{2}-k, 2^*k-1, k\gamma_q-1>0~\text{and}~ k\gamma_p-1<0.$$
We then obtain
\begin{align*}
M(a,b)+o_n(1)>&(\frac{1}{2}-k)(\|\nabla u_n\|_2^2+\|\nabla v_n\|_2^2)+\frac{\mu_1}{p}(k\gamma_p-1)\|u_n\|_p^p\\
\geq&(\frac{1}{2}-k)(\|\nabla u_n\|_2^2+\|\nabla v_n\|_2^2)+\frac{\mu_1}{p}(k\gamma_p-1)\|u_n\|_2^{p-\gamma_p}(\|\nabla u_n\|_2^2+\|\nabla v_n\|_2^2)^{\frac{\gamma_p}{2}}.
\end{align*}
Since $\frac{\gamma_p}{2}<1$, it follows that $\|\nabla u_n\|_2^2+\|\nabla v_n\|_2^2<+\infty$, establishing boundedness in $E$. Claim (i) follows directly.

By \eqref{eq:240729-1}, we have
$$
\lambda_{1,n}=\frac{1}{a}\left(-\|\nabla u_n\|_2^2+\mu_1\|u_n\|_p^p+\nu\alpha\int_{\R^N}u_n^{\alpha}v_n^{\beta}{\rm d}x\right)+o_n(1)
$$
and
$$
\lambda_{2,n}=\frac{1}{b}\left(-\|\nabla v_n\|_2^2+\mu_2\|v_n\|_q^q+\nu\beta\int_{\R^N}u_n^{\alpha}v_n^{\beta}{\rm d}x\right)+o_n(1).
$$
Since $\{(u_n,v_n)\}$ is bounded in $E$, it follows that $\{(\lambda_{1,n},\lambda_{2,n})\}$ is bounded, proving (ii). From (i) and (ii), claim (iii) follows by standard arguments.
\ep
\bl\label{le:240801-4}
Under the hypotheses of Corollary \ref{le:240801-2}. Let $(u_0,v_0)$ and $(\bar\lambda_1,\bar\lambda_2)$ be as in Lemma \ref{le:240801-1}. If additionally
\beq\label{eq:240801-12}
m_q^{\mu_2}(b)\geq \frac{2}{N-2}\nu^{-\frac{N-2}{2}}\alpha^{-\frac{(N-2)\alpha}{4}}\beta^{-\frac{(N-2)\beta}{4}}S^{\frac{N}{2}},
\eeq
then $(u_n,v_n)\to (u_0,v_0)$ strongly in $D_0^{1,2}(\R^N)\times D_0^{1,2}(\R^N)$. Furthermore, $u_0 \geq 0$ with $u_0 \not\equiv 0$, $v_0 \geq 0$ with $v_0 \not\equiv 0$, and $\bar\lambda_1, \bar\lambda_2 > 0$.
\el
\bp
We adopt the notations from Lemma \ref{le:240729-3} (see \eqref{eq:240801-5}-\eqref{eq:240827-1}) to facilitate concentration analysis. By \eqref{eq:240729-1} and the boundedness of $\{(\lambda_{1,n},\lambda_{2,n})\}$,
\beq\label{eq:240801-7}
\|\xi_1\|=\nu\alpha\|\eta_3\|, \quad \|\xi_2\|=\nu\beta\|\eta_3\|.
\eeq
Combining \eqref{eq:240801-9}, \eqref{eq:240827-1} and \eqref{eq:240801-7}, direct computation yields that if $\|\eta_3\| \neq 0$, then $\|\eta_3\|\geq\nu^{-\frac{N}{2}}\alpha^{\frac{(2-N)\alpha}{4}}\beta^{\frac{(2-N)\beta}{4}}S^{\frac{N}{2}}$.

Combining \eqref{eq:240801-3} and Concentration Compactness Principle, we obtain
\begin{align*}
P(u_n,v_n)=&\|\nabla u_0\|_2^2+\|\nabla v_0\|_2^2+\|\xi_1\|+\|\xi_2\|+\xi_{\infty,1}+\xi_{\infty,2}-\frac{\mu_1}{p}\gamma_p\|u_0\|_p^p-\frac{\mu_2}{q}\gamma_q\|v_0\|_q^q\\
&-2^*\nu\int_{\R^N}|u_0|^{\alpha}|v_0|^{\beta}{\rm d}x-2^*\nu\|\eta_3\|\\
=&P(u_0,v_0)+\|\xi_1\|+\|\xi_2\|+\xi_{\infty,1}+\xi_{\infty,2}-2^*\nu\|\eta_3\|\\
=&o_n(1).
\end{align*}
Since $P(u_0,v_0)=0$, it follows that
$$\|\xi_1\|+\|\xi_2\|+\xi_{\infty,1}+\xi_{\infty,2}=2^*\nu\|\eta_3\|.$$
Then combining \eqref{eq:240801-7}, we derive  $\xi_{\infty,1}=\xi_{\infty,2}=0$.

From \eqref{eq:240801-2}, we have
\beq\label{eq:240801-8}
\begin{aligned}
M(a,b)=&I(u_n,v_n)+o_n(1)\\
=&\frac{1}{2}(\|\nabla u_0\|_2^2+\|\nabla v_0\|_2^2+\|\xi_1\|+\|\xi_2\|)-\frac{\mu_1}{p}\|u_0\|_p^p-\frac{\mu_2}{q}\|v_0\|_q^q\\
&-\nu\left(\int_{\R^N}|u_0|^{\alpha}|v_0|^{\beta}{\rm d}x+\|\eta_3\|\right)+o_n(1)\\
=&I(u_0,v_0)+\frac{1}{2}(\|\xi_1\|+\|\xi_2\|)-\nu\|\eta_3\|+o_n(1)\\
=&I(u_0,v_0)+\frac{2^*-2}{2}\nu\|\eta_3\|+o_n(1).
\end{aligned}
\eeq

\textbf{Claim:} $\|\eta_3\| = 0$. Lemma \ref{le:240921-2} establishes that $m(a,b)=C_{a,b}=\inf\limits_{\cp_{a,b}}I(u,v)$ decreases in $a>0$ and $b\geq0$. If $\|\eta_3\| \neq 0$, then \eqref{eq:240801-8} and Lemma \ref{le:240801-1}-(iii) yield:
\begin{align*}
M(a,b)=&I(u_n,v_n)+o_n(1)=I(u_0,v_0)+\frac{2^*-2}{2}\nu\|\eta_3\|+o_n(1)\\
\geq&m(\|u_0\|_2^2,\|v_0\|_2^2)+\frac{2^*-2}{2}\nu\|\eta_3\|+o_n(1)\\
\geq&m(a,b)+\frac{2}{N-2}\nu^{-\frac{N-2}{2}}\alpha^{-\frac{(N-2)\alpha}{4}}\beta^{-\frac{(N-2)\beta}{4}}S^{\frac{N}{2}}+o_n(1),
\end{align*}
contradicting \eqref{eq:240801-10}. Thus $\|\xi_1\| = \|\xi_2\| = \xi_{\infty,1} = \xi_{\infty,2} = 0$, implying $(u_n,v_n) \to (u_0,v_0)$ in $D^{1,2}_0(\mathbb{R}^N)\times D^{1,2}_0(\mathbb{R}^N)$.

\textbf{Case 1: $u_0 \equiv 0$, $v_0 \equiv 0$.} Then $I(u_0,v_0)=0$. So \eqref{eq:240801-8} gives $M(a,b)=\frac{2^*-2}{2}\nu\|\eta_3\|=0$, contradicting $M(a,b) \geq k_0 > 0$.

\textbf{Case 2: $u_0 \equiv 0$, $v_0 \not\equiv 0$.} By Lemma \ref{le:240727-2},
$$I(u_0,v_0)=I(0,v_0)=m_q^{\mu_2}(\|v_0\|_2^2)\geq m_q^{\mu_2}(b)>0.$$
Then by \eqref{eq:240801-12} and \eqref{eq:240801-8}, it follows that
$$
M(a,b)\geq m_q^{\mu_2}(b)\geq \frac{2}{N-2}\nu^{-\frac{N-2}{2}}\alpha^{-\frac{(N-2)\alpha}{4}}\beta^{-\frac{(N-2)\beta}{4}}S^{\frac{N}{2}},
$$
contradicting \eqref{eq:240801-10}.

\textbf{Case 3: $u_0 \not\equiv 0$, $v_0 \equiv 0$.} By Lemma \ref{le:240727-2} and \eqref{eq:240801-8}:
$$M(a,b)=I(u_0,0)=m_p^{\mu_1}(\|u_0\|_2^2)<0,$$
contradicting $M(a,b)\geq k_0>0$.

Thus $u_0\gneqq0,v_0\gneqq0$. Finally, since $(u_0,v_0)$ solves \eqref{eq:240722-1}, Liouville's theorem \cite[Lemma A.2]{Ikoma2014} gives $\bar{\lambda}_1 > 0$, $\bar{\lambda}_2 > 0$.
\ep

\noindent
\textbf{Proof of Theorem \ref{th:240722-2}.}
By Lemma \ref{le:240727-2}, there exists $b_0>0$ such that for any $b < b_0$, inequality \eqref{eq:240801-12} holds. For any such $b < b_0$ and $\nu > 0$, the definition of $T_{a,b}$ implies that there exists a sufficiently small $a(b,\nu) > 0$ such that $T_{a,b} < c_0$ whenever $a < a(b,\nu)$. Consequently, Theorem \ref{th:240722-1} guarantees the existence of a normalized ground state solution $(u,v)$. Furthermore, Lemma \ref{le:20240927-1} provides the existence of a nonnegative Palais-Smale-Pohozaev sequence $\{(u_n, v_n)\}$.

Lemma \ref{le:240801-4} combined with \eqref{eq:240801-12} proves that the sequence $(u_n, v_n)$ converges strongly to $(u_0, v_0)$ in $D_0^{1,2}(\R^N) \times D_0^{1,2}(\R^N)$ and that $\bar\lambda_1 > 0, \bar\lambda_2 > 0$. Combining \eqref{eq:240729-1} and \eqref{eq:240731-1}, we obtain $\bar\lambda_1(a - \|u_0\|_2^2) = 0$ and $\bar\lambda_2(b - \|v_0\|_2^2) = 0$. This implies $\|u_0\|_2^2 = a$ and $\|v_0\|_2^2 = b$. Thus, the sequence $(u_n, v_n)$ also converges to $(u_0, v_0)$ in $E$. Finally, Corollary \ref{le:240801-2} yields
$$
I(u,v) < 0 < I(u_0, v_0) < I(u,v) + \frac{2}{N-2}\nu^{-\frac{N-2}{2}} \alpha^{-\frac{(N-2)\alpha}{4}} \beta^{-\frac{(N-2)\beta}{4}} S^{\frac{N}{2}}.
$$
\hfill$\Box$
		

\end{document}